\documentclass[9pt,twocolumn,twoside,fleqn]{extarticle}
\pdfoutput=1
\usepackage{mathpazo} 
\usepackage[scaled]{helvet}
\usepackage[T1]{fontenc}
\usepackage[left=48pt,right=42pt,top=46pt,bottom=60pt,headheight=15pt,headsep=10pt,letterpaper,twoside]{geometry}
\usepackage{graphicx}
\usepackage{authblk}
\usepackage{amsmath}
\usepackage{amssymb}
\usepackage{amsthm}
\usepackage{algorithm}
\usepackage{algorithmic}
\usepackage{dsfont}
\usepackage{colortbl}
\usepackage{xcolor}
\usepackage{booktabs}
\usepackage[protrusion=true,expansion=true]{microtype}
\usepackage{chngcntr}
\usepackage{hyperref}
\newfloat{algorithm}{t}{lop}

\restylefloat{algorithm}

\definecolor{color1}{RGB}{0,0,0} 
\definecolor{Gray}{gray}{0.95}
\definecolor{color2}{RGB}{33,64,154} 

\makeatletter
\newcommand\fs@coloruled{\def\@fs@cfont{\bfseries}\let\@fs@capt\floatc@coloruled
  \def\@fs@pre{\kern2pt}%
  \def\@fs@post{\kern2pt{\color{color1}\hrule}\relax}%
  \def\@fs@mid{\kern2pt{\color{color1}\hrule}\kern2pt}%
  \let\@fs@iftopcapt\iftrue}
\makeatother

\floatstyle{coloruled}

\newcommand{\realN}{\mathbb{R}}
\newcommand{\pinteN}{\mathbb{Z}_{+}}

\newcommand{\meter}{\text{m}}

\newcommand{\EX}{\mathbb{E}}
\newcommand{\PR}{\mathbb{P}}

\newcommand{\mT}{\text{T}}
\DeclareMathOperator*{\argmin}{argmin}

\newcommand{\vecOne}{\mathds{1}}
\DeclareMathOperator{\diag}{diag}
\DeclareMathOperator{\TV}{TV}

\newcommand{\trnYi}{Y_{\iota}^{\text{trn}}}

\newcommand{\tstYi}{Y_{\iota}^{\text{tst}}}

\newcommand{\nuAlgOne}{\hat{\nu}^{\,alg-1}}
\newcommand{\tauAlgOne}{\hat{\tau}^{\,alg-1}}
\newcommand{\betaAlgOne}{\hat{\beta}^{\,alg-1}}
\newcommand{\muAlgOne}{\hat{\mu}^{\,alg-1}}

\newcommand{\nuAlgTwo}{\hat{\nu}^{\,alg-2}}
\newcommand{\tauAlgTwo}{\hat{\tau}^{\,alg-2}}
\newcommand{\betaAlgTwo}{\hat{\beta}^{\,alg-2}}
\newcommand{\muAlgTwo}{\hat{\mu}^{\,alg-2}}

\newcommand{\tauAlgThree}{\hat{\tau}^{\,alg-3}}
\newcommand{\betaAlgThree}{\hat{\beta}^{\,alg-3}}
\newcommand{\muAlgThree}{\hat{\mu}^{\,alg-3}}

\title{A New Approach to Inverting Backscatter and Extinction from Photon-Limited Lidar Observations}

\author[1,2,$\dagger$]{Willem~J. Marais}
\author[2]{Robert~E. Holz}
\author[1]{Yu~Hen~Hu}
\author[2]{Ralph~E. Kuehn}
\author[3]{Edwin~E. Eloranta}
\author[1]{Rebecca~M. Willett}

\affil[1]{Department of Electrical and Computer Engineering (ECE), University of Wisconsin-Madison, Madison, Wisconsin, USA}
\affil[2]{Cooperative Institute for Meteorological Satellite Studies (CIMSS), University of Wisconsin-Madison, Madison, Wisconsin, USA}
\affil[3]{Space Science and Engineering Center (SSEC), University of Wisconsin-Madison, Madison, Wisconsin, USA}
\affil[$\dagger$]{Corresponding author: willem.marais@ssec.wisc.edu}

\begin{document}

\twocolumn[
\maketitle
\begin{@twocolumnfalse}
\begin{abstract}
\normalsize Atmospheric lidar observations provide a unique capability to directly observe the vertical column of cloud and aerosol scattering properties. Detector and solar background noise, however, hinder the ability of lidar systems to provide reliable backscatter and extinction cross-section estimates. Standard methods for solving this inverse problem are most effective with high signal-to-noise ratio observations that are only available at low-resolution in uniform scenes. This paper describes a novel method for solving the inverse problem with high-resolution, lower signal-to-noise ratio observations that are effective in non-uniform scenes. The novelty is twofold. First, the inference of the backscatter and extinction are done on images, whereas current lidar algorithms only use the information content of single profiles. Hence, the latent spatial and temporal information in the noisy images are utilized to infer the cross-sections. Second, the noise associated with photon-counting lidar observations can be modeled using a Poisson distribution, and state-of-the-art tools for solving Poisson inverse problems are adapted to the atmospheric lidar problem. It is demonstrated through photon-counting High Spectral Resolution Lidar (HSRL) simulations that the proposed algorithm yield inverted backscatter and extinction cross-sections (per unit volume) with smaller Mean Squared Error (MSE) values at higher spatial and temporal resolutions, compared to the standard approach. Two case studies of real experimental data are also provided where the proposed algorithm is applied on HSRL observations and the inverted backscatter and extinction cross-sections are compared against the standard approach.
\end{abstract}
\setcounter{tocdepth}{1}
\tableofcontents
\vspace{7mm}
\end{@twocolumnfalse}
]

\section{Introduction}
\label{sec:introduction}

Atmospheric lidar systems measure the profiles of attenuated backscatter. From these the profiles of backscatter and extinction cross-sections are inferred. Examples of lidar systems are the NASA-CNES Cloud-Aerosol Lidar and Infrared Pathfinder Satellite Observations (CALIPSO) mission, the NASA Langley High Spectral Resolution Lidar (HSRL) instrument, the NASA Cloud Physics Lidar (CPL) instrument, the Micro-Pulsed Lidar (MPL) network and the upcoming ESA-JAXA ATmospheric LIDar (ATLID) which will be on the Earthcare satellite \cite{winker2009overview,hair2008airborne,durand2007esa,mcgill2002cloud,welton2001global}. Another example of ground-based lidar is the Wisconsin-Madison HSRL system, which is the first system for which we have developed scatter inference tools, with the intention to extend it to other lidar instruments in the near future. All these lidar systems provide the earth science community with the unique capability to resolve the atmospheric vertical structure with very high sensitivity.

For the most part the lidar signal of interest is typically weak when observing tenuous aerosols and clouds. The implication is that the photon detectors employed are typically photon starved, and the detectors produce noisy energy rate measurements accompanied with spurious detections due to the dark current \cite{liu2006estimating}. Furthermore, the solar background is another noise source and it further increases the difficulty to detect small changes in the signal. An extra complication is that the lidar detector noise is in most cases non-Gaussian, and the noise variance is signal dependent \cite{liu2006estimating,liu2002simulation}. Thus, statistical estimators well known in the atmospheric science community, such as the ``optimal estimator'' or Tikhonov regularization\footnote{Tikhonov regularization is also known as ridge regression.} \cite{hastie2009elements,rodgers2000inverse}, cannot be directly applied unless a substantial amount of averaging is employed so to invoke the central limit theorem with a strong set of assumptions.

The predominant methodology used by lidar experts to estimate the cross-sections is to 1) divide the noisy images into non-overlapping blocks, 2) average the noisy observations in each block to reduce the noise variance, 3) solve for the unknowns using a single scatter lidar equation for each averaged block and~4) employ 1-D lowpass filters across the temporal and range axes after the calculations to suppress residual noise \cite{eloranta2014high,nishizawa2008algorithm,young2008caliopATBDp4,ansmann2007particle}. The size of each block is not always fixed, it depends on the amount of averaging that is required to reduce the variance (uncertainty) of the estimates to a satisfactory level. Averaging can introduces unnecessary biases, since most lidar observations of clouds and aerosols are not spatially uniform. Lowpass filtering can sometimes yield satisfactory results, but it is not suited for features that have high frequency components, such as clouds that have sharp boundaries. Furthermore, current parameter inversion techniques are ill-equipped to reliably invert the desired physical parameters. For example, in standard backscatter lidar systems such as CALIOP a type of fixed point iteration algorithm is used which is not always numerically stable \cite[\S~2.3.2.3]{young2008caliopATBDp4}.

\subsection{Contributions}

The primary contribution of this paper is the introduction of statistical estimators that use the noise physical model in conjunction with the single scatter lidar equation, to fit the estimates of the backscatter and extinction cross-section \emph{images} on the noisy observations. Estimates of the backscatter and extinction cross-section images are constrained to be spatially piecewise smooth. In other words, the relation between the cross-section at a specific pixel and the neighboring cross-section pixels, is controlled by a constraint that limits the variation among the cross-sections. This happens while the new approach simultaneously finds the optimal solution for a large domain of observations, using both the noise physical model and non-linear lidar model. This contribution is novel, since current lidar algorithms only consider a single profile when inferring the backscatter and extinction cross-sections. Although block averaging does impose a smoothness constraint when it is used to reduce the noise variance, its smoothness constraint is per block and not among the blocks wherein the average quantities are calculated.

With new approach we present in this paper 1) the spatial and temporal information of the noisy observations is utilized and 2) the discontinuities and high frequency components of the backscatter and extinction cross-sections are preserved. The algorithms presented in this paper are based on sound theory that guarantees uniqueness of an estimate under reasonable conditions. We demonstrate that with the new approach the extinction and backscatter cross-sections can be inferred more accurately compared to the standard approach, with both simulated and real data. We specifically focus on inverting the particulate backscatter and extinction cross-sections from low SNR photon-counting HSRL observations. The noise of photon-counting systems can be accurately modeled by a Poisson distribution \cite{liu2006estimating}. Our intention is to extend our research to standard backscatter lidar systems, and other photon detection systems in the near future, such as analog mode PMT lidar systems, where the noise is Compounded Poisson \cite{liu2002simulation}.

The statistical estimator that we employ is known as the Total Variation (TV) Penalized Maximum Likelihood Estimator (PMLE). Compelling results have been obtained with the TV-PMLE in medical imaging systems \cite{oh2013logTV,harmany2012spiral}. Some medical imaging systems, such as Positron Emission Tomography (PET), low-dose Computed Tomography (CT) and confocal microscopy, also use photon-counting detectors and have the same noise characteristics as photon-counting lidar system \cite{willett2003platelets}. We demonstrate, through simulated observations and two case studies on real data, that the backscatter and extinction cross-sections can be inverted at a lower SNR with smaller Mean Squared Error (MSE) and normalized MSE values compared to the standard inversion approach. We achieve this by using the algorithm Sparse Poisson Image Reconstruction ALgorithm (SPIRAL), and an adaptation of the algorithm log-SPIRAL \cite{oh2013logTV,harmany2012spiral}. These algorithms are well established and are numerically stable under reasonable assumptions, such as the convexity of the optimization problems. 

Our secondary contribution is a novel approach to estimating the lidar ratio\footnote{The lidar ratio is also known as the extinction-to-backscatter ratio \cite[\S~2.3.1.1]{young2008caliopATBDp4}.}, compared to the customary approach of estimating the extinction cross-section directly. The approach we propose has the benefit of exploiting the inverted backscatter cross-section to constrain the estimation of the lidar ratio and the extinction cross-section. 

Although we mainly focus on ground based HSRL photon-counting systems, we believe that the research we present is a necessary step to further improve upon the backscatter and extinction inversion techniques for all types of standard backscatter lidar systems, such as the space-based analog lidar systems CALIOP. In particular, it is necessary to consider the noise model if further improvements are to be made on backscatter and extinction inversion techniques. Furthermore, since lidar instruments produce images, the backscatter and extinction cross-sections should be inferred using all the available spatial and temporal information in the images. 

\subsection{Paper outline}

The paper starts off in \S\ref{sec:prelimary} with an introduction to the HSRL models. Next in \S\ref{sec:standard_approach} we describe the standard approach to infer the unknown coefficients and other contributions that had been made in inferring the unknown parameters. Our proposed approach is presented in \S\ref{sec:new_approach}. In \S\ref{sec:experimental_results} we present two simulation results to compare the performances of the standard and new approaches against each other. Thereafter in \S\ref{sec:exprmnts_actual_data} we present two case studies to demonstrate the capabilities of the proposed approach. The paper ends with the conclusion and discussion of future work in \S\ref{sec:conclussion}.

\subsection{Assumptions}
\label{sec:intro:assumptions}

We primarily focus on the parallel (polarized) backscatter cross-section estimated from the parallel polarized photon-counts. It is assumed that an estimate of the linear depolarization coefficient is available, which is used to compute the backscatter cross-section from the parallel backscatter cross-section. The assumed linear depolarization coefficient is denoted by $\rho\in[0,1]^{N\times{}K}$; see equation 2.3 of~\cite{weitkamp2006lidar} for a definition of $\rho$.

We assume that multiple-scattering is negligible. Without this assumption the inference problem becomes extremely complex and we defer to future research to take in account multiple scattering with the proposed new approach. For now we assume that the solid angle of the lidar receiver is small enough so that multiple scattering is less significant, e.g. 100 micro-radians for a ground-based HSRL~\cite{eloranta2014high}~\cite[Chapter~5]{weitkamp2006lidar}.

It is assumed that the lidar ratio is greater or equal to one throughout the text. Examples can be constructed using orientated platelet ice-crystals to show that this assumption is not always true, especially when the angle of the laser axis is perpendicular to the earth's surface. If the angle is oblique, then the assumption is reasonable.

\subsection{Notation convention and symbols}
\label{subsec:notation_convention}

A raw lidar image consists of $N$ range bins (row axis) and $K$ profiles (column axis); a profile refers to a single column. The range and profiles are expressed in meters and seconds. The row axis of all processed images, i.e. extinction cross-section estimates, is converted to altitude of mean sea level (msl). A 15km by 2~hour lidar image, with a resolution of 7.5m by 2.5s, has 2000 range bins and 2880 profiles.

The set of non-negative real numbers is denoted by $\realN_{+}$, and the set of non-negative integers is denoted by $\pinteN$. To save space, we squeeze a double summation term $\sum_{n=1}^{N}\sum_{k=1}^{K}$ into one $\sum_{n=1,k=1}^{N,K}$. The vector $e_{n}\in\realN^{N}$ is a canonical vector. The symbol $\vecOne_{N}\in\realN^{N}$ represents a vector of $N$ ones. 

In several of the models described in this paper, it is more convenient to use pointwise multiplication operations than linear algebra multiplication operations. The symbol $\cdot$ indicates that two matrices are multiplied pointwise. For example, if $A,B\in\realN^{N\times{}K}$, we have that 
\begin{equation}
    e_{n}^{\mT}[A\cdot{}B]e_{k} = [e_{n}^{\mT}Ae_{k}][e_{n}^{\mT}Be_{k}],
\end{equation}
where with a slight abuse of notation the $e_{n}$ and $e_{k}$ vectors are of different lengths. 

We denote the particulate parallel backscatter and extinction cross-sections (per unit volume) by the symbols $\nu\in\realN_{+}^{N\times{}K}$ and $\beta\in\realN_{+}^{N\times{}K}$; these are the primary unknown parameters of interest which are introduced in \S\ref{sec:prelimary}. The backscatter cross-section is denoted by the symbol $\nu_{+}$, where the subscript~$+$ indicates that the backscatter cross-section is the sum of the parallel and perpendicular backscatter cross-sections. The backscatter cross-section is also known as the backscatter volume coefficient in the realm of atmospheric science; the same applies to the extinction cross-section \cite{young2008caliopATBDp4,petty2006first}. Whenever $\beta$ is estimated with say algorithm number one, its estimate will be denoted by $\betaAlgOne$; the same applies to $\nu$.

The matrix $Q\in\realN^{N\times{}N}$ is a lower triangular matrix of ones, scaled by a constant $\Delta{}r$. The constant $\Delta{}r$ is the range sampling resolution of the lidar instrument. Each row of $Q$ represent the Riemann integral:
\begin{align}
    e_{n}^{\mT}[Q\beta{}]e_{k} 
    & = \Delta{}r\sum_{l=1}^{n}e_{l}^{\mT}\beta{}e_{k}.
\end{align}
Hence, when $Q$ acts on a matrix $\beta$, the output is the scaled cumulative sum of each column of $\beta$.

The Poisson noisy observations of the HSRL molecular- and combined-channels (defined in \S\ref{subsec:hsrl_phy_model}) are denoted by $Y_{m}\in\pinteN^{N\times{}K}$ and $Y_{c}\in\pinteN^{N\times{}K}$. The letter $m$ is an abbreviation of the word molecular, and $c$ is an abbreviation of the word combined. The symbol $\iota\in\{c, m\}$ is used to index $c$ and $m$. Each entry in $Y_{\iota}$ is statistically independent from each other.

To aid the reader interpret the text and the equations, Table~\ref{tbl:symbol_table} gives a non-exhaustive list of symbols that are used in paper. A short description of each symbol is given, the type of variable (matrix or scalar) and the location where the symbol was introduced.
\begin{table}[!ht]
    \caption{This table gives a non-exhaustive list of symbols that are used in paper. A short description of each symbol is given, the type of variable (matrix or scalar) and the location where the symbol was introduced.}
    \label{tbl:symbol_table}
    \begin{tabular}{@{}p{0.7cm}>{\raggedright\arraybackslash}p{4.2cm}cl>{\raggedright\arraybackslash}p{0.5cm}@{}} \toprule
        Symbol    & Description                                                  & Type & Location \\ \midrule
        $\nu$     & Parallel backscatter cross-section                           & Matrix & \S\ref{sec:prelimary} \\
        $\beta$   & Extinction cross-section                                     & Matrix & \S\ref{sec:prelimary} \\
        $\mu$     & Lidar ratio                                                  & Matrix & \S\ref{sec:prelimary} \\
        $\nu_{+}$ & Non-polarized backscatter cross-section                      & Matrix & \S\ref{sec:prelimary} \\
        $\tau$    & Optical depth                                                & Matrix & \S\ref{sec:prelimary} \\
        $\rho$    & Depolarization coefficient                                   & Matrix & \S\ref{sec:intro:assumptions} \\
        $S_{c},S_{m}$ & Combined- and molecular-channel HSRL models & Matrix     & \S\ref{subsec:hsrl_phy_model} \\
        $Y_{c},Y_{m}$ & Combined- and molecular-channel HSRL photon-count images & Matrix & \S\ref{subsec:notation_convention}, \S\ref{subsec:poisson_noise_model}  \\
        $C_{g}$   & Gain calibration parameter                                   & Matrix & \S\ref{subsec:hsrl_phy_model} \\
        $C_{mc}$  & Molecular backscatter calibration parameter of combined-channel & Matrix & \S\ref{subsec:hsrl_phy_model} \\
        $C_{am}$  & Particulate backscatter calibration parameter of molecular-channel & Scalar & \S\ref{subsec:hsrl_phy_model} \\
        $C_{mm}$  & Molecular backscatter calibration parameter of molecular-channel & Matrix  & \S\ref{subsec:hsrl_phy_model} \\
        $b_{c}$   & Dark and solarbackground counts of combined-channel & Matrix & \S\ref{subsec:hsrl_phy_model} \\
        $b_{m}$   & Dark and solarbackground counts of molecular-channel & Matrix & \S\ref{subsec:hsrl_phy_model} \\
        $Q$       & Integrator matrix & Matrix & \S\ref{subsec:notation_convention} \\\arrayrulecolor{color2}\bottomrule
    \end{tabular}
\end{table}

\section{Preliminary - The unknowns \& the HSRL and noise models}
\label{sec:prelimary}

The unknown parameters that we consider in this paper are the particulate extinction $\beta\in\realN_{+}^{N\times{}K}$ and parallel (polarized) backscatter $\nu\in\realN_{+}^{N\times{}K}$ cross-sections; the unit of both these coefficients is $\meter^{-1}$. The backscatter cross-section is computed using the linear depolarization coefficient
\begin{equation}
    \nu_{+} \equiv \nu / (\vecOne - \rho),
\end{equation}
and it is assumed that we already have an estimate of the depolarization measurements; the division is taken to be pointwise.

In addition to the unknowns $\beta$ and $\nu$, we are also interested in the optical depth and the lidar ratio. The optical depth is denoted by the symbol $\tau\in\realN_{+}^{N\times{}K}$ and its relation to $\beta$ is $\tau \equiv Q\beta$, where $Q\in\realN^{N\times{}N}$ represents integration (see \S\ref{subsec:notation_convention}). The lidar ratio, the ratio between the extinction and backscatter cross-sections, is denoted by 
\begin{equation}
    \mu \equiv \beta / \nu_{+}.
\end{equation}

In the following subsection visual examples of the unknowns are given. In the next subsection the lidar models for the High Spectral Resolution Lidar (HSRL) system are introduced. HSRL models, which are based on the single scatter lidar equation, are typically written as continous functions that are indexed mainly by a range or an altitude index variable~\cite{hair2008airborne,weitkamp2006lidar,eloranta2005high}. We'll deviate from this convention, since the new approach that we are presenting works with images and also require full knowledge of all the calibration parameters. After the HSRL models are introduced, the noise physical model is introduced.

\subsection{Non-uniformity and smoothness properties of the unknown parameters}

\begin{figure}[!ht]
    \centering
    \includegraphics[width=\linewidth]{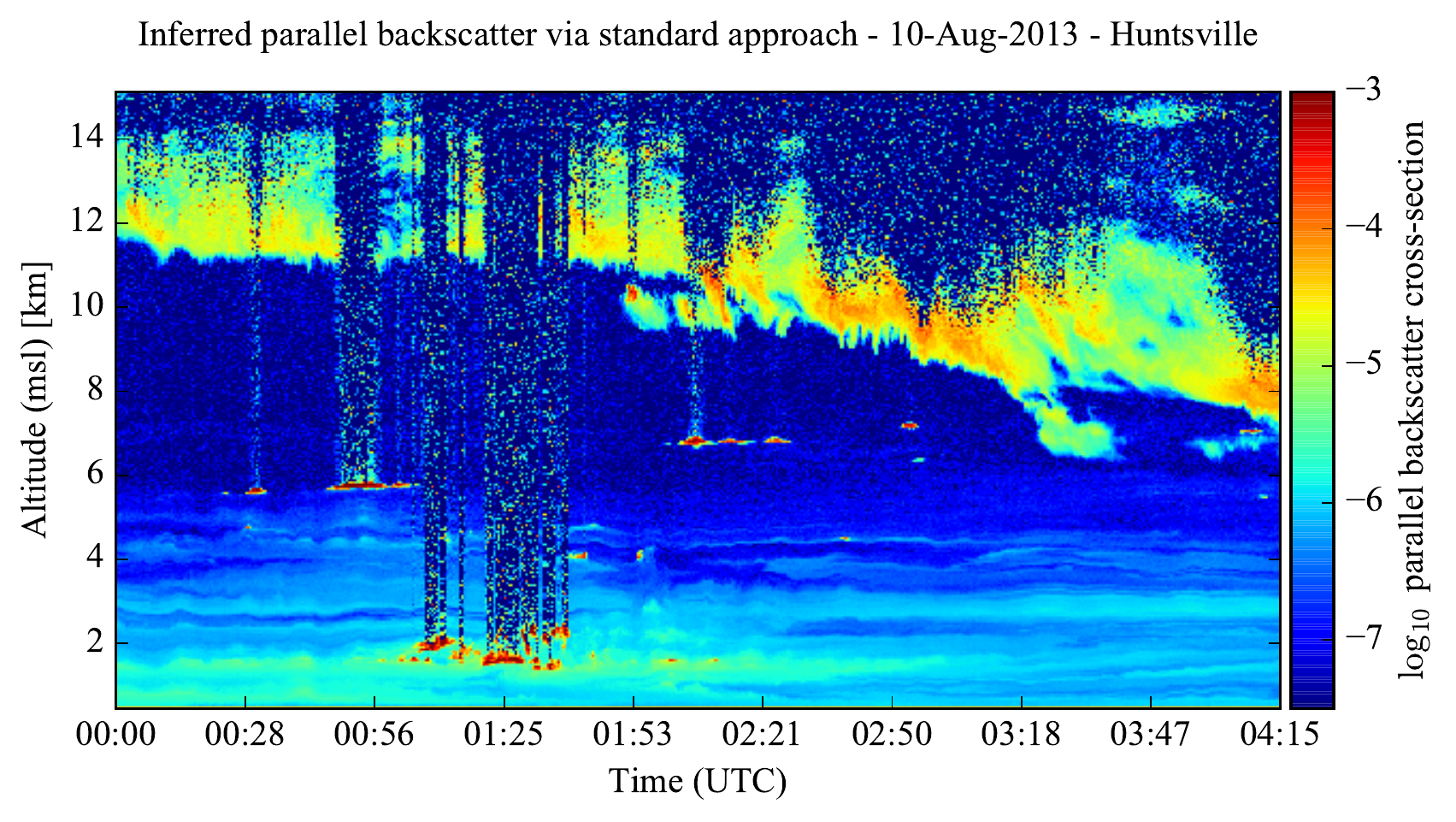}
    \caption{An example of the inverted particulate parallel backscatter cross-section of the Wisconsin-Madison HRSL system \cite{eloranta2014high,hsrldev-www}. A significant amount of averaging was done to increase the SNR. Each hour consists of 770 columns and each 1~km range consists of about 133 rows. The vertical noisy stripes are due to the complete attenuation of the laser pulse by the clouds.}
    \label{fig:noisy_backscatter_image}
\end{figure}

Figure~\ref{fig:noisy_backscatter_image} shows an example image of inverted particulate parallel backscatter $\nu$ using the standard approach algorithm (see \S\ref{sec:standard_approach}); the lidar observations are from the Wisconsin-Madison ground based HSRL instrument \cite{eloranta2014high,hsrldev-www} (see \S\ref{subsec:hsrl_phy_model}). This purpose of Figure~\ref{fig:noisy_backscatter_image} is to show the non-uniformity of a typical scene. The HSRL instrument is stationary, therefore we see the backscatter of the atmosphere as it moved across the instrument.
The backscattered energy image was averaged from a resolution of 7.5m (row-axis) by 2.5s (column-axis) to 30m by 30s, thus the detector noise is not clearly visible. Between 0km and 3km at 1:25~UTC various small dense clouds are present, which have backscatter cross-section of about $10^{-3}$ or larger. To the left and right of these low level clouds are faint aerosol layers which have backscatter cross-section ranging from $10^{-5}$ to $10^{-6}$. Various ripples are present on top of the aerosol layers and various non-uniform horizontal layers are present below 4km. Above 8km a large cloud is present that has backscatter cross-section that ranges from $10^{-6}$ to $10^{-3}$.

\begin{figure}[!ht]
    \centering
    \includegraphics[width=\linewidth]{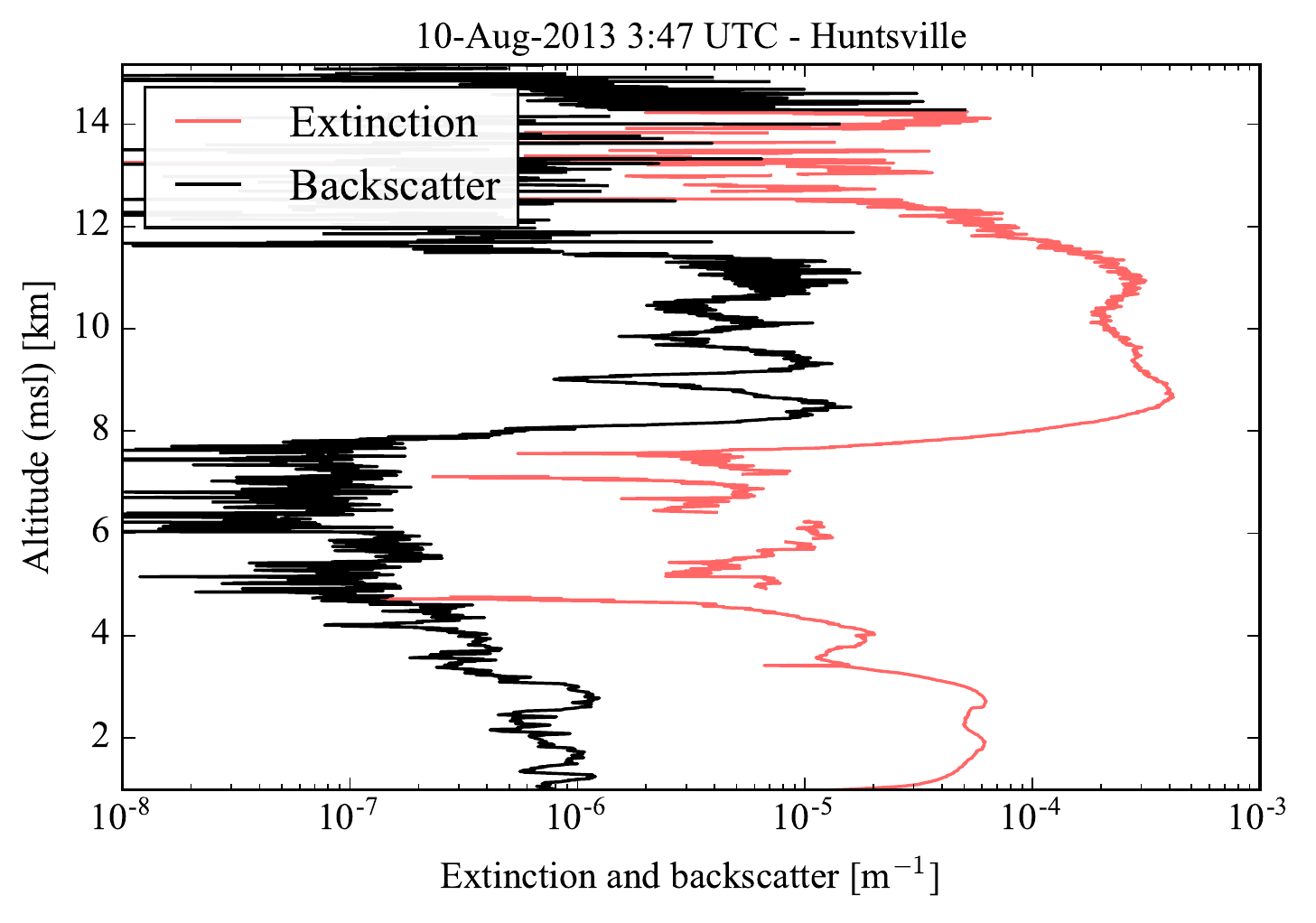}
    \caption{The black graph shows the inverted backscatter $\nu$ cross-section of the 1min averaged HSRL observations, corresponding to Figure~\ref{fig:noisy_backscatter_image}, using the standard approach algorithm (see \S\ref{sec:standard_approach}). The red graph shows the corresponding inverted extinction $\beta$ cross-section.}
    \label{fig:noisy_profile_backscatter_and_scatter}
\end{figure}

Figure~\ref{fig:noisy_profile_backscatter_and_scatter} show examples of the inverted backscatter and extinction cross-sections of a single profile, using the standard approach algorithm (see \S\ref{subsec:hsrl_phy_model} and \S\ref{sec:standard_approach}). The lidar observations were averaged over 1min to reduce the noise in the estimates. 
The inverted extinction cross-section is much smoother than the backscatter cross-sections, since the standard approach use a lowpass filter that is applied on the altitude axis to reduce residual noise. Figure~\ref{fig:noisy_profile_backscatter_and_scatter} gives an idea of how the backscatter and extinction cross-sections are related. The extinction cross-section can be about two orders of magnitude larger than the backscatter cross-sections, and it more or less have the same upward and downward trends compared to the backscatter cross-section. From the graphs in Figure~\ref{fig:noisy_profile_backscatter_and_scatter} we can deduce that the underlying backscatter and extinction cross-sections are piecewise smooth functions.

\subsection{The HSRL models}
\label{subsec:hsrl_phy_model}

An HSRL that we consider in this paper has three channels, which we will refer to as the combined, molecular and cross-polarization-channels~\cite{eloranta2005high}. With the combined-channel the parallel (polarized) backscattered energy of both particulates and molecules are measured, and with the molecular-channel primarily the molecular backscattered energy is measured. This is achieved by using an iodine filter as a notch-filter, to reject the particulate backscatter so that primarily the molecular backscatter Doppler broadened signal is measured~\cite{hair2008airborne,eloranta2005high}\cite[p. 262]{petty2006first}. The cross-polarization-channel measures the perpendicular backscattered energy of both particulates and molecules. With these three channels, the backscatter and extinction cross-sections along with the depolarization coefficient can be estimated. The HSRL cross-polarization-channel model will not be discussed, since it is assumed that an estimate of the linear depolarization coefficient is available; refer to \cite{weitkamp2006lidar} for more information.

The combined- and molecular-channel models derived from the single scatter lidar equation~\cite{weitkamp2006lidar}. The combined-channel model, which maps the parallel backscatter $\nu$ and extinction $\beta$ cross-section matrices to a backscattered energy image, is defined by
\begin{equation}
    S_{c}(\nu,\beta) = C_{g}\cdot(\nu + C_{mc})\cdot\exp\left(-2Q\beta\right) + b_{c} \label{eq:physical_model_comb_channel}.
\end{equation}
The only unknowns in this model are the parallel backscatter~$\nu$ and extinction~$\beta$ cross-section matrices, and the rest are precomputed calibration matrices. The symbol $\cdot$ indicates that the matrices are multiplied pointwise; why this is necessary will become clear in the next paragraph. The matrix $C_{g}\in\realN_{+}^{N\times{}K}$ is the gain calibration matrix which includes the transmitted laser energy, the receiver solid angle which is a function of the telescope area, the geometric overlap function, the molecular transmittance and the optical-system and detector efficiency coefficients~\cite[Chapter~5]{weitkamp2006lidar}; refer to~\cite{eloranta2014high} for more information about how the geometric overlap function is estimated. The parallel molecular backscatter cross-section and the calibration parameter that modifies it, is represented by $C_{mc}\in\realN_{+}^{N\times{}K}$. The mentioned calibration parameter describes the portion of the molecular backscatter that is attenuated by the solarbackground bandpass filter. The parallel polarized solarbackground energy is denote by $b_{c}\in\realN_{+}^{N\times{}K}$. Note that the columns of~$b_{c}$ change as a function of the temporal axis, since the solarbackground radiation change as a function of time. 

The reason why pointwise and not linear algebra matrix multiplication is used in the HSRL model, is due to the calibration matrix $C_{g}$. Each column of $C_{g}$ represents a different profile, and there are at least two calibration parameters that can change as a function of the profile number. For a ground-based lidar system, such as the Wisconsin-Madison HSRL instrument, the photon detectors are saturated for a duration after the laserpulse has been transmitted~\cite{razenkov2010characterization}. Thus, a calibration parameter is the particulate transmittance between the time at which the laserpulse was transmitted and the first altitude bin at which the backscattered energy is measured. Since the first altitude bin of $S_{c}(\nu,\beta)$ is when the photon detectors are not saturated, the calibration matrix $C_{g}$ includes the unobserved particulate transmittance. For space-based lidar systems where a photon counting image can span a large geographical area, the molecular transmittance can change as a function of the profile number in $C_{g}$. Therefore, the HSRL model as it is written in \eqref{eq:physical_model_comb_channel} models a lidar image accurately.

The HSRL molecular-channel model is defined by 
\begin{equation}
    S_{m}(\nu,\beta) = C_{g}\cdot(C_{am}\nu + C_{mm})\cdot\exp\left(-2Q\beta\right) + b_{m} \label{eq:physical_model_molec_channel}.
\end{equation}
The molecular-channel model is defined so that its gain calibration matrix is equal to that of the combined-channel, which is~$C_{g}$. The scalar $C_{am}\in\realN_{+}$ represents the rejection of the particulate backscatter by the iodine filter. The calibration matrix $C_{mm}\in\realN_{+}^{N\times{}K}$ represents parallel molecular backscatter cross-section and the calibration parameter that modifies it. The mentioned calibration parameter describes the portion of the molecular backscatter that is attenuated by the solarbackground bandpass filter and the iodine filter. The solarbackground energy is denote by $b_{m}\in\realN_{+}^{N\times{}K}$.

If $\nu$ and $\beta$ are the true parameters of an HSRL scene, we have that $\EX[Y_{c}] = S_{c}(\nu,\beta)$ and $\EX[Y_{m}] = S_{m}(\nu,\beta)$. To simplify the mathematical expressions of the HSRL models in the rest of this paper, we defined the matrices $S_{c}$ and $S_{m}$ to be implicit functions of the relavant unknowns parameters:
\begin{equation}
    S_{c} \equiv S_{c}(\nu,\beta)\,\text{ and }\,S_{m} \equiv S_{m}(\nu,\beta).
\end{equation}

\begin{figure}[!ht]
    \centering
    \includegraphics[width=\linewidth]{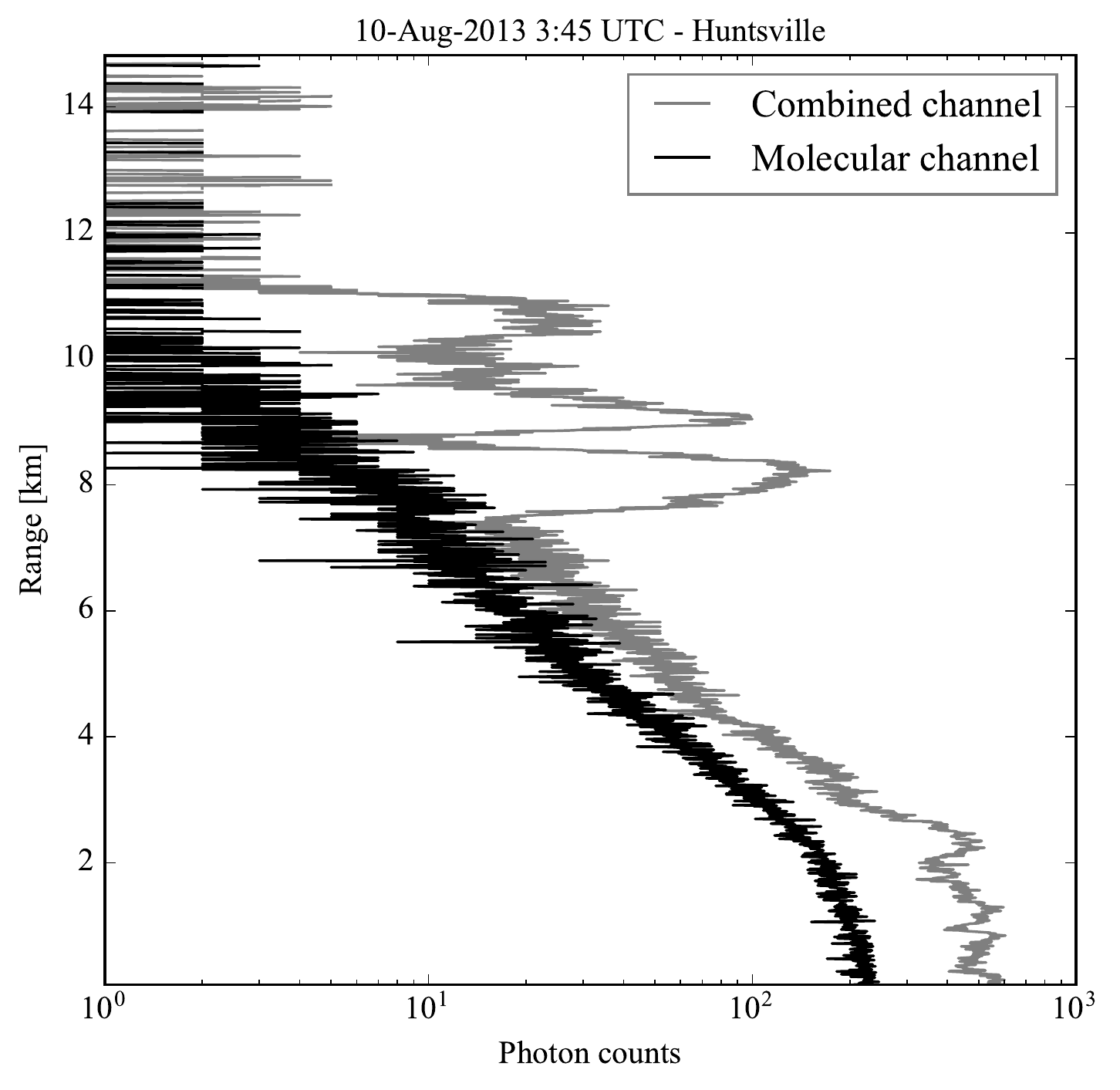}
    \caption{An example of noisy observations of the molecular $S_{m}$ and combined $S_{c}$ channels. These noisy observations are of a single column of the observations that were used to create Figures~\ref{fig:noisy_backscatter_image} and~\ref{fig:noisy_profile_backscatter_and_scatter}. The range resolution is~7.5m and the temporal resolution of the profile is~2.5s.}
    \label{fig:noisy_hsrl_observations}
\end{figure}

Figure~\ref{fig:noisy_hsrl_observations} shows an example of noisy observations of~$S_{m}$ and~$S_{c}$.
These noisy observations are of a single column of the observations that were used to create Figures~\ref{fig:noisy_backscatter_image} and~\ref{fig:noisy_profile_backscatter_and_scatter}. Observe that molecular-channel photon counts profile is smoother than combined-channel photon counts profile. This is because the molecular-channel primarily measures the attenuated backscatter of molecules, whereas the combined channel is very sensitive to the particulate backscatter. This is due to $C_{am} < 10^{-3}$. Also notice also that rate at which the photon counts decrease, which is due to the receiver solid angle and the attenuation of the laser pulse.

\subsection{The Poisson noise physical model}
\label{subsec:poisson_noise_model}

It is well-known that the noise of a photon-counting system can be modeled by a Poisson Probability Mass Function (PMF)~\cite{liu2006estimating,willett2003platelets}. Let $Y_{\iota}\in\pinteN^{N\times{}K}$ be the Poisson noisy observations of $S_{\iota}$, where $\iota\in\{c, m\}$. The Poisson PMF of $Y_{\iota}$ is defined as 
\begin{align}\label{eq:poisson_pmf}
    & \PR(Y_{\iota} = y) 
    = \prod_{n=1}^{N}\prod_{k=1}^{K}\exp(-S_{\iota,n,k})\frac{(S_{\iota,n,k})^{y_{n,k}}}{y_{n,k}!} \\
    & S_{\iota,n,k} = e_{n}^{\mT}\EX[Y_{\iota}]e_{k} \\
    & y_{n,k} = e_{n}^{\mT}(y)e_{k}, 
\end{align}
where $e_{n}$ and $e_{k}$ are canonical vectors and with a slight abuse of notation these two vector differ in lengths (see \S\ref{subsec:notation_convention}). The Penalized Maximum Likelihood Estimator (PMLE) is a well-known technique which uses the noise model to fit the desired parameters to the noisy observations, where the desired parameters are constrained to enforce some a priori information \cite{willett2007multiscale,willett2003platelets,eggermont2001maximum}. We show in \S\ref{sec:new_approach} how the HSRL models are use in conjunction with the Poisson noise model to infer the extinction $\beta$ and backscatter $\nu$ cross-sections.

\section{The standard approach in inverting the backscatter and extinction cross-sections}
\label{sec:standard_approach}

Ideally when there is no detector noise and the calibration parameters are fully known, the extinction $\beta$ and parallel backscatter $\nu$ cross-sections can be solved algebraically~\cite{eloranta2014high}. From the HSRL models~\eqref{eq:physical_model_comb_channel} and \eqref{eq:physical_model_molec_channel} we have that
\begin{align}\label{eq:standard_retrieve_beta_p}
    \tau = Q\beta 
         = 
    -\frac{1}{2}\log_{e}\left[
    \frac{(S_{c} - b_{c})C_{am} - (S_{m} - b_{m})}
    {C_{g}\cdot(C_{mc}C_{am} - C_{mm})}
    \right]
\end{align}
and
\begin{equation}\label{eq:standard_retrieve_parl_nu_p}
    \nu =
    \frac
    {(S_{c} - b_{c})\cdot{}C_{mm} - (S_{m} - b_{m})\cdot{}C_{mc}}
    {S_{m} - b_{m} - (S_{c} - b_{c})C_{am}}.
\end{equation}
The extinction cross-section is computed by applying a discrete derivative operator on the optical depth $\tau$. 

Algorithm~\ref{alg:standard_approach_algorithm} concisely describes the necessary steps to infer the backscatter and extinction cross-section, together with the optical depth and the lidar ratio.
\begin{algorithm}
\caption{The standard approach in inferring the backscatter $\nu_{+}$ and extinction $\beta$ cross-sections, the lidar ratio $\mu$ and optical depth $\tau$. The tuning parameters of $F^{\text{Avg}}$ and $F^{\text{Low}}$ are set by a lidar expert or an automated system, which determined by some objective criteria. Note that whenever the unknown parameters are estimated with this algorithm, it is indicated by a superscript $alg-1$.}
\label{alg:standard_approach_algorithm}
    Let $F^{\text{Low}}(\cdot)$ be a lowpass filter function, and let $F^{\text{Avg}}(\cdot)$ be an averaging function. And let $F^{\text{Diff}}(\cdot)$ be a discrete derivative operator.
\begin{algorithmic}[1]
    \STATE $S_{c}^{\text{Avg}} = F^{\text{Avg}}(Y_{c})$
    \STATE $S_{m}^{\text{Avg}} = F^{\text{Avg}}(Y_{m})$
    \STATE $\hat{\tau}^{\,alg-1} \leftarrow $ plug $S_{c}^{\text{Avg}}$ and $S_{m}^{\text{Avg}}$ into \eqref{eq:standard_retrieve_beta_p} 
    \STATE $\nuAlgOne \leftarrow$ plug $S_{c}^{\text{Avg}}$ and $S_{m}^{\text{Avg}}$ into \eqref{eq:standard_retrieve_parl_nu_p}
    \STATE $\betaAlgOne = F^{\text{Diff}}(F^{\text{Low}}(\hat{\tau}^{\,alg-1}))$
    \STATE $\hat{\nu}_{+}^{\,alg-1} = \nuAlgOne / (\vecOne - \rho)$
    \STATE $\muAlgOne = (\vecOne - \rho) \cdot \betaAlgOne / \hat{\nu}^{\,alg-1}$
    \RETURN $\hat{\nu}_{+}^{\,alg-1}, \betaAlgOne, \muAlgOne, \tauAlgOne$
\end{algorithmic}
\end{algorithm}
Figure~\ref{fig:noisy_profile_backscatter_and_scatter} shows an example of inverted parallel backscatter and extinction cross-sections of Algorithm~\ref{alg:standard_approach_algorithm}.

Between lines 1 and 2 the noise variance of the measured backscattered energy images is reduced using an averaging function $F^{\text{Avg}}(\cdot)$. This averaging function either divides an image into non-overlapping blocks and calculate the average statistic for each block, or a moving average filter is used which is equivalent to a two-dimensional convolution of the image with a normalized two-dimensional rectangular function. The type of average function that is used, depends uniformity of the features in the scene. The amount of averaging that is required, depends on what the required standard deviation (uncertainty) of the estimates should be. Between lines~3 and~5 estimates of the optical depth and parallel backscatter cross-section are computed. At line~5 an estimate of the extinction cross-section is computed, by first applying a low-pass filter function $F^{\text{Low}}(\cdot)$ on the estimated optical depth and then a discrete derivative operator $F^{\text{Diff}}(\cdot)$ is applied. An example of low-pass filter function, is a Savitzky-Golay (SG) filter that is first applied on the temporal axis and then another SG filter on the range axis~\cite{eloranta2014high,schafer2010frequency}. The backscatter cross-section is computed at line~6, and then the lidar ratio is computed at line~7. The symbol $\vecOne$ represents a $N\times{}K$ matrix of ones. 

\subsection{An overview of other inversion methods}
\label{subsec:over_current_retrieve_meths}

Past research, similar to ours, has investigated the inference of particulate backscatter $\nu$ and extinction $\beta$ cross-sections. Our contribution relies on an HSRL system that has the benefit of not having to make a priori assumptions about the lidar ratio, whereas some of the cited contributions only use a standard backscatter lidar system. With a standard backscatter lidar system stronger assumptions have to be made about the relation between the extinction and backscatter cross-sections, since two unknowns have to be estimated from one set of observations. Hence, a direct comparison cannot be made between this paper's contribution and some of the contributions cited in this section. Nevertheless, we list the past contributions to highlight the advantages of the proposed method.

Averaging is the standard noise reduction method for lidar observations in combination of other noise suppressing techniques \cite{eloranta2014high,nishizawa2008algorithm,young2008caliopATBDp4,ansmann2007particle}. In \cite{eloranta2014high} a Savitzky-Golay lowpass filter is used suppress residual noise prior to retrieving the particulate extinction cross-section from HSRL observations (see step~3 of Figure~\ref{alg:standard_approach_algorithm}). This has the drawback of smoothing the desired signal itself and inadvertently introducing significant estimation biases to the inverted cross-sections on top of the bias introduced by averaging. The algorithm developed in \cite{young2008caliopATBDp4} tries to avoid biasing due to averaging, by selectively applying different amounts of averaging on pre-classified features. In our contribution we try to do the least amount of accumulations of photon counts, and use a penalty function which does not impose over-smoothing constraints on the estimated cross-sections and introduce less bias. We demonstrate in \S\ref{sec:experimental_results} and \S\ref{sec:exprmnts_actual_data} the backscatter can be inverted at a higer resolution than the standard approach.

In \cite{pornsawad2012retrieval} the physical model was transformed into linear system which acts on a extinction cross-section vector. The extinction cross-section vector is solved for by employing the iterative Levenberg-Marquardt algorithm. The solution for the extinction vector is regulated by choosing the number of iterations of the algorithm, and this is achieved by using the L-curve method which assumes that the noise is Gaussian \cite{pornsawad2012retrieval,hansen1999curve}. Although it is not explicitly stated, it seems that it is implicitly assumed in \cite{pornsawad2012retrieval} that the noise is Gaussian distributed whereas in our approach we use a more accurate Poisson noise model.

A recent suggested approach relies on a transform method to remove high frequency components from the measured energy rates of a standard backscatter lidar system \cite{tian2014improved}. The transform method, referred to as the empirical mode decomposition, decomposes the measured energy rates into low and high frequency components \cite{tian2014improved}. It appears that in \cite{tian2014improved} it is assumed that the noise can have negative values and therefore is not Poisson, even though the lidar instrument that \cite{tian2014improved} refers to employs photon-counting \cite{cimel2015ce370-2}.

\section{A new approach to infer the backscatter and extinction cross-sections}
\label{sec:new_approach}

Our goal with the new approach is to invert the extinction $\beta$ and backscatter $\nu$ cross-sections at lower SNRs compared to the standard approach. In other words a smaller number of photon-counts have to be accumulated, vertically and temporally, to increase the SNR of the observations in order to invert for $\beta$ and~$\nu$. We achieve this by ensuring that the estimated cross-sections adhere to several constraints. A Total Variation (TV) smoothness constraint is imposed on the parallel backscatter~$\nu$ and extinction $\beta$ cross-sections, which enforce the estimates to be spatially piecewise constant. With block averaging the estimate of a block is calculated independently from the rest of the blocks, whereas with the TV constraint the whole image is taken in account. The TV constraint ensures that the estimates are piecewise constant with few discontinuities. It imposes no restrictions on where the discontinuities are and hence allows high-resolution estimates.

In addition to imposing a smoothness constraint on $\nu$ and $\beta$, we want to use the inverted parallel backscatter cross-sections~$\nu$ along with the depolarization coefficient $\rho$ to constrain the inversion of the extinction cross-sections $\beta$. We motivate this approach by noting that in the presence of Poisson noise, the estimation of $\beta$ is statistically ill-posed \cite{o1986statistical}. In other words, when the estimation of $\beta$ is not constrained, any small changes in the observations - due to for example Poisson noise - can lead to significant changes to an estimate of $\beta$.

We achieve our goals by employing a well known technique: the Penalized Maximum Likelihood Estimator (PMLE). The basic idea of the PMLE is that we seek estimates which 1) are a good fit to the observed data as measured using the Poisson likelihood in \eqref{eq:poisson_pmf} and 2) adhere to our a priori model of piecewise smoothness. This idea can be formulated as an optimization problem in which we search over all candidate cross-sections, and choose the cross-sections which minimize the sum of the negative Poisson log likelihood and a penalty term which is smaller for piecewise smooth cross-sections with small TV. The MLE use the noise model to fit the desired parameters on the noisy observations and it is based on well established theorems~\cite{eggermont2001maximum}. The TV seminorm has been used in several other signal-processing applications to promote piecewise smoothness with satisfactory results \cite{oh2013logTV,needell2013stable,harmany2012spiral}. A regularization parameter is used to set the degree to which the smoothness is promoted, and the regularization parameter is automatically set using a cross-validation heuristic \cite{harmany2012spiral}. We adopt the algorithms SPIRAL and log-SPIRAL to compute the PMLE  \cite{oh2013logTV,harmany2012spiral}.

\subsection{Motivation for using the Total Variation (TV) smoothness constraint}
\label{subsec:motivate_TV}

\begin{figure}[!ht]
    \centering
    \includegraphics[width=\linewidth]{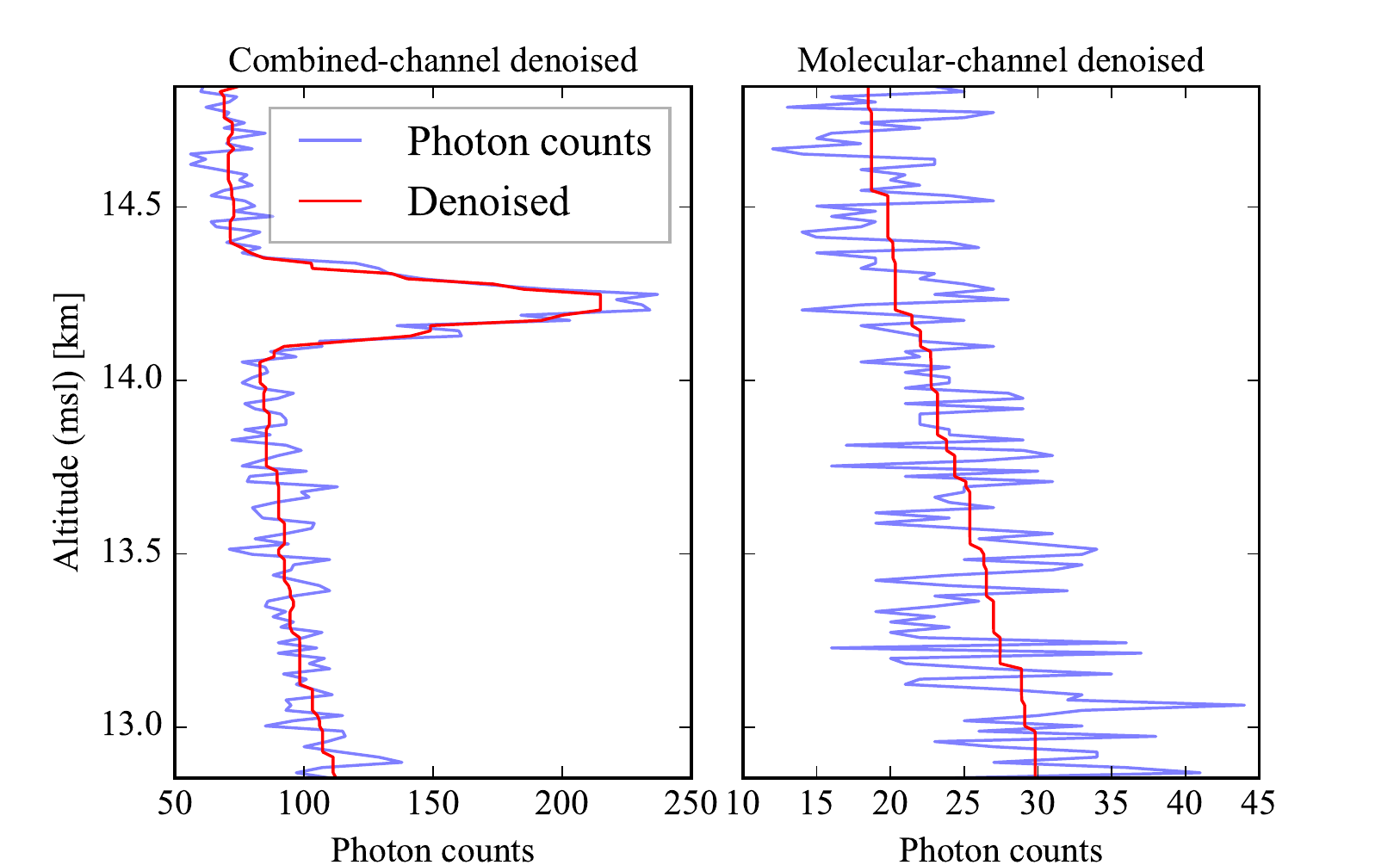}
    \caption{The images in this figure show what is meant by approximating a piecewise smooth signal by a piecewise constant function, when denoising the Poisson noisy observations of the corresponding piecewise smooth signal. The light blue graphs in the figure show the combined- and molecular-channel Poisson noisy photon-counts of the Wisconsin-Madison HSRL~\cite{eloranta2014high}. Between 14~km and 14.5~km a cirrus cloud is present. The red graphs show the denoised backscattered energy, estimated by the SPIRAL-TV algorithm~\cite{harmany2012spiral}. From the left image we see that the SPIRAL-TV is able to find the appropriate locations of the cirrus cloud discontinuities. The right image shows how the underlying molecular signal is appropriated by a piecewise constant function.}
    \label{fig:denoise_example}
\end{figure}

Figures~\ref{fig:noisy_backscatter_image} and~\ref{fig:noisy_profile_backscatter_and_scatter} suggest that the unknown parallel backscatter and extinction cross-section are spatially and temporally piecewise smooth. When TV is used as a smoothness constraint, the piecewise constant function that it induces can yield accurate and robust approximation results of noisy spatial piecewise smooth functions, for both Poisson and Gaussian noisy images~\cite{needell2013stable,harmany2012spiral,chan2005image}. Figure~\ref{fig:denoise_example} shows what is meant by approximating a piecewise smooth signal by a piecewise constant function, when denoising the Poisson noisy observations of the corresponding piecewise smooth signal.
The light blue graphs in the figure show the combined- and molecular-channel Poisson noisy photon-counts of the Wisconsin-Madison HSRL~\cite{eloranta2014high}. Between 14~km and 14.5~km a cirrus cloud is present. The red graphs show the denoised backscattered energy, estimated by the SPIRAL-TV algorithm~\cite{harmany2012spiral}. From the combined-channel denoised backscattered energy, we can see that the SPIRAL-TV algorithm was able to detect where the discontinuities are below and above the cirrus cloud. And from the molecular-channel denoised backscattered energy, the underlying monotonically decreasing molecular signal is approximated by a piecewise constant function, similar to a decreasing staircase function. The length of each step is automatically determined by the TV-PMLE.

\subsection{Formulation of TV-PMLE optimization problem and the algorithm}
\label{subsec:formulate_tv_pmle}

First in \S\ref{subsubsec:reparameterization_nu} we describe the formulation of the TV-PMLE when estimating the parallel backscatter cross-section $\nu$. The inverted backscatter is then used to constrain the estimation of the extinction cross-section $\beta$, which is described in the following subsection \S\ref{subsubsec:reparameterization_beta}.

\subsubsection{Inversion of the backscatter cross-section}
\label{subsubsec:reparameterization_nu}

To estimate the parallel backscatter cross-section, the HSRL combined- and molecular-channels are denoised and a modified version of \eqref{eq:standard_retrieve_parl_nu_p} is used to compute the parallel backscatter cross-section. For each channel we estimate the backscattered energy using the background energy as a priori information. Let $\omega_{\iota}\in\realN_{+}^{N\times{}K}$ represent the backscattered energy in the absence of the background energy for channel $\iota\in\{c,m\}$, where $c$ and $m$ are abbreviations of the words combined and molecular. In other words, we have that 
\begin{equation}
    \omega_{\iota} \equiv S_{\iota}(\nu,\beta) - b_{\iota};
\end{equation}
see \eqref{eq:physical_model_comb_channel} and \eqref{eq:physical_model_molec_channel}. An example of $\omega_{c}$ and $\omega_{m}$ estimates, are the red graphs of Figure~\ref{fig:denoise_example} without the background energies of the combined- and molecular-channels. We can now reparameterize the backscattered energy model $S_{\iota}(\nu,\beta)$ with the model 
\begin{equation}
    f_{\iota}(\omega_{\iota}) = \omega_{\iota} + b_{\iota}.\label{eq:backscattered_eng_model}
\end{equation}
Once we have an estimates of $\omega_{c}$ and $\omega_{m}$, we can use an algebraic expression similar to \eqref{eq:standard_retrieve_parl_nu_p} to compute an estimate of the parallel backscatter cross-section~$\nu$; we have that 
\begin{equation}
    \nu = 
    \frac{\omega_{c}\cdot{}C_{mm} - \omega_{m}\cdot{}C_{mc}}
    {\omega_{m} - \omega_{c}C_{am}}.
\end{equation}
Once we have an estimate of the parallel backscatter cross-section, we can compute the backscatter cross-section $\nu_{+}$ using the depolarization coefficient $\rho$: $\nu_{+} = \nu / (\vecOne - \rho)$.

With the MLE technique the estimate $\hat{\omega}_{\iota}$ is obtained, for a given $Y_{\iota}$, by minimizing the negative log-likelihood of \eqref{eq:poisson_pmf} with respect to the constraint $\omega_{\iota}\in\realN_{+}^{N\times{}K}$. A penalty function is added to the MLE to regularize the estimate of $\omega_{\iota}$. If it is reasonable to assume that $\omega_{\iota}$ can be approximated with piecewise constant functions, then an adequate penalty function is the TV seminorm. The specific TV seminorm that is used in this paper is
\begin{equation}\label{eq:tv_semi_norm}
    \|\omega\|_{\text{TV}}
    = \sum_{\substack{n=1,\\k=1}}^{\substack{N-1,\\K}} \left|\omega_{n,k} - \omega_{n+1,k}\right|
    + \sum_{\substack{n=1,\\k=1}}^{\substack{N,\\K-1}} \left|\omega_{n,k} - \omega_{n,k+1}\right|,
\end{equation}
where the subscripts in this case are indices to the rows and columns of the matrix $\omega$. This variant of the TV seminorm is known as the anisotropic TV seminorm.

The TV-PMLE is formulated by
\begin{align}
    & \hat{\omega}_{\iota} = \argmin_{\omega\in\Omega_{\iota}}
        \left\{\ell(\omega; Y_{\iota}) + \lambda\|\omega\|_{\text{TV}}\right\} \label{eq:opt_objective_function} \\
    & \ell(\omega; Y_{\iota})
        = \vecOne_{N}^{\mT}\left[f_{\iota}(\omega) - Y_{\iota}\cdot\log_{e}f_{\iota}(\omega)\right]\vecOne_{K}, \label{eq:mle_loss_function}
\end{align}
where $\lambda \geq 0$ is referred to as the regularizer parameter. The vectors $\vecOne_{N}\in\realN^{N}$ and $\vecOne_{K}\in\realN^{K}$ are vectors of ones, and the symbol $\cdot$ indicates that the matrices are multiplied pointwise. Throughout this paper \eqref{eq:mle_loss_function} is referred to as the loss function. A cross-validation heuristic is used to choose an appropriate regularizer parameter~$\lambda$~\cite{oh2013logTV}. The cross-validation technique divides the photon counting observations into two sets of noisy observations using Poisson thinning (see \S\ref{app:poisson_thinning}). With the first set of noisy observations the regularizer parameter $\lambda$ is varied, and for each $\lambda$ an estimate of $\hat{\omega}_{\iota}$ is obtained using \eqref{eq:opt_objective_function}. Then with the second set of noisy observations the reconstructed backscattered energy (i.e. plugging $\hat{\omega}_{\iota}$ into~\eqref{eq:backscattered_eng_model}) is validated; see Algorithm~\ref{alg:cross-validation} and~\cite{oh2013logTV} for more information.

The minimizer $\hat{\omega}_{\iota}$ is unique if the negative log-likelihood is strictly convex, since the TV seminorm is convex. This will hold true if none of the elements of $Y_{\iota}$ are zero. The algorithms SPIRAL and log-SPIRAL guarantee convergence to a solution if~1) the gradient of the loss function \eqref{eq:mle_loss_function} is Lipschitz continuous and~2) the objective function \eqref{eq:opt_objective_function} is coercive \cite[\S{}III.D]{harmany2012spiral}. The loss function \eqref{eq:mle_loss_function} satisfies both these conditions; refer to \cite{harmany2012spiral} for more details. 

\subsubsection{Inversion of the extinction cross-section}
\label{subsubsec:reparameterization_beta}

The inverted parallel backscatter cross-section can be used to constrain the inversion of the extinction cross-section. We do this by estimating the lidar ratio using the parallel backscatter cross-section, since there is a linear relation between the lidar ratio and the extinction cross-section. Once we obtain an estimate of the extinction cross-section, we can then algebraically compute the estimate of the optical depth.

The HSRL has both parallel and perpendicular polarized channels, which measure the backscattered energy with separate detectors. The relative gain of the two channels have to be reliably calibrated to get the depolarization coefficient in order to obtain the non-polarized backscatter cross-section. For this application, the inversion procedure for the extinction cross-section is applied on the parallel channels and inferred parameters are then corrected for the depolarization dependence.

To estimate the lidar ratio, we reparameterize the molecular HSRL model~\eqref{eq:physical_model_comb_channel}. Let $g(\cdot)$ be the reparameterization of $S_{m}(\nu,\beta)$, which is defined as 
\begin{align}
    & g(\tilde{\mu}) = C_{bm}\cdot\exp\left(-2Q\left[\hat{\nu}\cdot\tilde{\mu}\right]\right) + b_{m} \label{eq:model_g_mu_hat_nu} \\
    & \mu = \tilde{\mu}\cdot(\vecOne - \rho) \label{eq:mu_pseudo_mu} \\
    & C_{bm} = C_{g}\cdot(\hat{\nu}C_{am} + C_{mm}) \label{eq:matrix_G_iota},
\end{align}
where $\tilde{\mu}\in[1, \infty)^{N\times{}K}$. \eqref{eq:mu_pseudo_mu} shows where the estimate of the lidar ratio is corrected for the depolarization dependence. The MLE technique is used to estimate the lidar ratio by minimizing the negative log-likelihood of \eqref{eq:poisson_pmf}. The TV seminorm is used a penalty function to regularize the estimate of $\tilde{\mu}$. The TV-MPLE is formulated by 
\begin{align}
    & \hat{\tilde{\mu}} = \argmin_{\tilde{\mu}\in[1, \infty)^{N\times{}K}}
        \left\{\ell(\tilde{\mu}; Y_{m}) + \lambda\|\tilde{\mu}\|_{\text{TV}}\right\} \label{eq:mu_opt_objective_function} \\
    & \ell(\tilde{\mu}; Y_{m})
        = \vecOne_{N}^{\mT}\left[g(\tilde{\mu}) - Y_{m}\cdot\log_{e}g(\tilde{\mu})\right]\vecOne_{K}, \label{eq:mu_mle_loss_function}
\end{align}
where the loss function $\ell(\cdot; \cdot)$ has been redefined and $\lambda \geq 0$ is the regularizer parameter. The loss function \eqref{eq:mu_mle_loss_function} is not necessarily strictly convex for all $\tilde{\mu}\in[1, \infty)^{N\times{}K}$, unless $b_{m} = 0$. Subsection \S\ref{subsubsec:uniqueness_scatter_coeff} describes for what values of $\tilde{\mu}$ is \eqref{eq:mu_mle_loss_function} strictly convex. Also the gradient of \eqref{eq:mu_mle_loss_function} is not Lipschitz continuous, since $\exp(\cdot)$ is not a Lipschitz continuous function. Thus $\tilde{\mu}$ has to be constrained to a closed subset of $[1, \infty)^{N\times{}K}$ so that the conditions can be met under which the log-SPIRAL algorithm converge to a minimizer. 

\subsection{Discussion about the loss function}

Since the noise we are dealing with is Poisson, the loss function we use is different from what is typically used in optimal estimation or Tikhonov regularization~\cite{rodgers2000inverse}. A Gaussian Probability Density Function (PDF) could be used to model the noise. When we are dealing with low photon counts, a Poisson PMF arguably describes the noise better than a Gaussian PDF. Hence, it is more accurate to use the actual noise model when deriving the loss function to do the estimation of the unknown cross-sections.

\subsection{The algorithm and details regarding it}
\label{subsec:proposed_algorithm}

A succinct outline of the new algorithm is given in Algorithm~\ref{alg:new_approach_nu_beta_mu_tau}. 
\begin{algorithm}
\caption{The new approach in inverting for the backscatter $\nu_{+}$ and extinction $\beta$ cross-sections, along with the lidar ratio $\mu$ and the optical depth $\tau$. Note that whenever the unknown parameters are estimated with this algorithm, it is indicated by a superscript $alg-2$.}
\label{alg:new_approach_nu_beta_mu_tau}
\begin{algorithmic}[1]
    \STATE \COMMENT{Invert the parallel backscatter cross-section}
    \STATE Set $\ell(\omega; Y_{\iota}) = \vecOne_{N}^{\mT}\left[f_{\iota}(\omega) - Y_{\iota}\cdot\log_{e}f_{\iota}(\omega)\right]\vecOne_{K}$
    \FORALL{$\iota\in\{c, m\}$}
        \STATE $\hat{\omega}_{\iota} = $ Algorithm~\ref{alg:cross-validation} [$\ell(\omega; Y_{\iota})$, $Y_{c}$, $Y_{m}$, $\realN_{+}^{N\times{}K}$]
    \ENDFOR
    \STATE $\nuAlgTwo = (\omega_{c}\cdot{}C_{mm} - \omega_{m}\cdot{}C_{mc})/(\omega_{m} - \omega_{c}C_{am})$
    \STATE $\hat{\nu}_{+}^{\,alg-2} = \nuAlgTwo / (\vecOne - \rho)$
    \STATE \COMMENT{Invert the extinction cross-section}
    \STATE Set $\ell(\tilde{\mu}; Y_{\iota}) = \vecOne_{N}^{\mT}\left[g(\tilde{\mu}) - Y_{\iota}\cdot\log_{e}g(\tilde{\mu})\right]\vecOne_{K}$
    \STATE $\hat{\tilde{\mu}} = $ Algorithm~\ref{alg:cross-validation} [$\ell(\tilde{\mu}; Y_{m})$, $Y_{m}$, $\mathcal{M}$]
    \STATE $\muAlgTwo = \hat{\tilde{\mu}}\cdot(\vecOne - \rho)$
    \STATE $\betaAlgTwo = \hat{\nu}_{+}^{\,alg-2}\cdot\muAlgTwo$
    \STATE $\tauAlgTwo = Q\betaAlgTwo$
    \RETURN $\hat{\nu}_{+}^{\,alg-2}, \betaAlgTwo, \muAlgTwo, \tauAlgTwo$
\end{algorithmic}
\end{algorithm}
To make Algorithm~\ref{alg:new_approach_nu_beta_mu_tau} more readable, the details regarding TV-PMLE and the cross-validation are presented in Algorithm~\ref{alg:cross-validation}. The parallel backscatter cross-section is inverted between lines~2 and~6 in Algorithm~\ref{alg:new_approach_nu_beta_mu_tau}, and the inverted parallel backscatter is denoted by $\nuAlgTwo$. The backscatter cross-section $\nu_{+}$ is computed at line~7, and is denoted by $\hat{\nu}_{+}^{\,alg-2}$. Then $\nuAlgTwo$ is used to invert for the lidar ratio between lines 9 and 11, and the inverted lidar ratio is denoted by $\muAlgTwo$. Thereafter the extinction cross-section and the optical depth are computed, and these are denoted by $\betaAlgTwo$ and $\tauAlgTwo$. 

\newcounter{mycount}
\counterwithin{algorithm}{mycount}
\refstepcounter{mycount}
\renewcommand\thealgorithm{\arabic{algorithm}-\alph{mycount}}
\setcounter{algorithm}{1}

\begin{algorithm}
\caption{The TV-PMLE with cross-validation to estimate either $\omega_{\iota}$ or $\mu$; see Algorithm~\ref{alg:new_approach_nu_beta_mu_tau}.}
\label{alg:cross-validation}
\begin{algorithmic}[1]
    \REQUIRE A loss function such as $\ell(\psi; Y_{c}, Y_{m})$, the Poisson noisy matrices $Y_{c}$, $Y_{m}$ and the constraint set $\Psi$ for $\psi$.
    \STATE $Y_{c}^{\text{trn}}, Y_{c}^{\text{tst}} \leftarrow \text{Poisson Thinning}\,(Y_{c})$
    \STATE $Y_{m}^{\text{trn}}, Y_{m}^{\text{tst}} \leftarrow \text{Poisson Thinning}\,(Y_{m})$
    \FORALL{$\lambda$ in $\Lambda$}
        \STATE $\displaystyle\hat{\psi}(\lambda) = \argmin_{\psi\in\Psi}\{\ell(\psi; Y_{c}^{\text{trn}}, Y_{m}^{\text{trn}}) + \lambda\|\psi\|_{\TV}\}$
    \ENDFOR
    \STATE $\displaystyle\lambda^{*} = \argmin_{\lambda\in\Lambda}\ell(\hat{\psi}(\lambda); Y_{c}^{\text{tst}}, Y_{m}^{\text{tst}})$
    \RETURN $\hat{\psi}(\lambda^{*})$
\end{algorithmic}
\end{algorithm}

\subsubsection{Parallel backscatter cross-section inversion algorithm - lines~2 to~7}

Algorithm~\ref{alg:cross-validation} is used to estimate $\hat{\omega}_{c}$ and $\hat{\omega}_{m}$, in order to estimate the parallel backscatter cross-section. From lines 1 and 5, of Algorithm~\ref{alg:cross-validation}, a cross-validation heuristic is used to find the regularizer parameters that best constrain the estimates $\hat{\omega}_{\iota}$. The cross-validation heuristic requires two independent copies of the random matrix $Y_{\iota}$, and these can be created using Poisson thinning (see appendix~\S\ref{app:poisson_thinning}) \cite{oh2013logTV}. In lines~1 and~2 the Poisson noisy observations $Y_{\iota}$ are thinned into two matrices $Y_{\iota}^{\text{trn}}$ and $Y_{\iota}^{\text{tst}}$, where $Y_{\iota} = Y_{\iota}^{\text{trn}} + Y_{\iota}^{\text{tst}}$ and $\EX[Y_{\iota}^{\text{trn}}] = 0.5\times{}S_{\iota}(\nu,\beta)$. The superscripts trn and tst are abbreviations of the words train and test.

Between lines 3 and 5 the unknown $\hat{\omega}_{\iota}$ is estimated by varying the piecewise constant constraint with a regularizer parameter $\lambda$. The regularizer parameter $\lambda$ is selected from a list of non-negative numbers $\Lambda$; for example $\Lambda = \{10^{-2}, 10^{-1.8}, \ldots, 10^{0.8}, 10\}$. The optimization problem in line 4 is solved using the SPIRAL algorithm \cite{harmany2012spiral}.

In line 6 the regularization parameter is chosen that minimizes the loss function with $Y_{\iota}^{\text{tst}}$. Next, in line 6 of Algorithm~\ref{alg:new_approach_nu_beta_mu_tau} the parallel backscatter cross-section is computed using $\hat{\omega}_{c}$ and $\hat{\omega}_{m}$. And at line~7 the backscatter cross-section is computed using the parallel backscatter cross-section $\hat{\nu}^{\,alg-2}$ with the depolarization coefficient $\rho$.

\subsubsection{Extinction cross-section inversion algorithm - lines~8 to~13}

As with the parallel backscatter cross-section inversion algorithm, the lidar ratios are estimated by varying a regularizer parameter $\lambda$. The set $\Lambda$ will likely be different from what was used for the backscatter inversion algorithm, since lidar ratio smoothness constraint will be different. 

The loss function $\ell(\tilde{\mu}; Y_{m})$ is not necessarily strictly convex for all $\tilde{\mu}\in[1, \infty)^{N\times{}K}$, and $\tilde{\mu}$ has to be constrained to a closed set $\mathcal{M}\subset[1, \infty)^{N\times{}K}$. Furthermore the set $\mathcal{M}$ has to be bounded, since it is required that the gradient of the loss function is Lipschitz continuous. The next subsection \S\ref{subsubsec:uniqueness_scatter_coeff} defines the set $\mathcal{M}$. 

An adaptation of the log-SPIRAL algorithm is used to implement line 4 of Algorithm~\ref{alg:cross-validation} \cite{oh2013logTV}, where the gradient matrix of $\ell(\tilde{\mu}; Y_{\iota})$ is 
\begin{equation}
    \nabla\ell(\tilde{\mu}; Y_{\iota})
    = -2\hat{\nu}\cdot{}Q^{\mT}\left[(g(\tilde{\mu}) - b_{\iota})\cdot\left(\vecOne - \frac{Y_{\iota}}{g(\tilde{\mu})}\right)\right].
\end{equation}
Once the lidar ratio is estimated, in line~12 of Algorithm~\ref{alg:new_approach_nu_beta_mu_tau}, the extinction cross-section $\betaAlgTwo$ are estimated using both $\muAlgTwo$ and $\hat{\nu}_{+}^{\,alg-2}$. Thereafter the optical depth, $\tauAlgTwo$, is computed by cumulatively integrating every column of $\betaAlgTwo$.

\subsection{Ensuring uniqueness of the inverted extinction cross-section $\beta$}
\label{subsubsec:uniqueness_scatter_coeff}

We will now show that Algorithm~\ref{alg:new_approach_nu_beta_mu_tau} can only invert the extinction cross-section, with a guarantee that the estimate is unique, when the total optical depth less than a predefined threshold. Without loss of generality assume that $K = 1$. The Hessian matrix of the loss function $\ell(\tilde{\mu}; Y_{m})$ (line~9 of Algorithm~\ref{alg:new_approach_nu_beta_mu_tau}) is 
\begin{align}
    & H(\tilde{\mu}) = 4\diag\left\{\hat{\nu}\right\}Q^{\mT}D(\tilde{\mu})Q\diag\left\{\hat{\nu}\right\} \label{eq:hessian_matrix_mu} \\
    & D(\tilde{\mu}) = \diag\left\{(g(\tilde{\mu}) - b_{m})\cdot\left[\vecOne_{N} - \frac{Y_{m}\cdot{}b_{m}}{[g(\tilde{\mu})]^{2}}\right]\right\},
\end{align}
where the division is taken to be pointwise. The loss function $\ell(\tilde{\mu}; Y_{m})$ is strictly convex if the Hessian matrix $H$ is positive definite, and p.d. (positive definiteness) is ensured if 
\begin{align}
    & g(\tilde{\mu}) > \sqrt{Y_{m}\cdot{}b_{m}} \label{eq:func_g_lower_bound} \\
    \Rightarrow\, & \hat{\tau} = Q[\hat{\nu}\cdot\tilde{\mu}] < -\frac{1}{2}\log_{e}\left[\frac{\sqrt{Y_{m}\cdot{}b_{m}} - b_{m}}{C_{bm}}\right].\label{eq:tau_upper_bound}
\end{align}
In other words, there is an upper bound limit on the optical depth that can be uniquely estimated and the upper bound is influenced by the background energy rate $b_{m}$. This in turn restricts the range of extinction cross-section that can be uniquely estimated, for a specific lidar scene. The constraint set of $\tilde{\mu}$ is then 
\begin{equation}
    \mathcal{M} = \left\{\tilde{\mu}\in[1, \infty)^{N\times{}K} : H(\tilde{\mu}) \text{ is p.d.}\right\}.
\end{equation}
An easy, albeit crude, approach to determine $\mathcal{M}$, is to assume a single lidar ratio for the whole scene and then find the maximum lidar ratio for which the Hessian matrices are p.d. The disadvantage of this approach is that the lidar ratio upper bound can be too conservative, since clear-sky portions of a scene could have negligible lidar ratios. 

\subsection{Another algorithm worthy of notice to invert the extinction cross-section $\beta$}
\label{subsec:another_alg_est_beta}

\renewcommand\thealgorithm{\arabic{algorithm}}

It is not always clear whether it is better to regularize an estimate of the lidar ratio $\mu$, and then to compute the extinction cross-section~$\beta$. Or to regularize an estimate of $\beta$, and then compute the lidar ratio $\mu$. With this in mind, Algorithm~\ref{alg:different_approach_beta} shows another approach to invert the extinction cross-section $\beta$, where the reparameterized lidar model is 
\begin{equation}
    h(\beta) = C_{bc}\cdot\exp(-2Q\beta) + b_{m}.
\end{equation}
The matrix $C_{bc}$ is defined in \eqref{eq:matrix_G_iota}. In this case the extinction cross-section is constrained to be greater or equal to the inverted backscatter cross-section.
\begin{algorithm}
\caption{Another approach in inverting for the particulate extinction cross-section $\beta$. Note that whenever the unknown parameters are estimated with this algorithm, it is indicated by a superscript $alg-3$.}
\label{alg:different_approach_beta}
\begin{algorithmic}[1]
    \REQUIRE Previously inverted parallel backscatter cross-section $\hat{\nu}$, such as $\nuAlgTwo$, with the depolarization coefficient $\rho$ to compute $\hat{\nu}_{+}$.
    \STATE \COMMENT{Invert the extinction cross-section}
    \STATE Set $\ell(\beta; Y_{m}) = \vecOne_{N}^{\mT}\left[h(\beta) - Y_{\iota}\cdot\log_{e}h(\beta)\right]\vecOne_{K}$
    \STATE $\betaAlgThree = $Algorithm~\ref{alg:cross-validation} [$\ell(\beta; Y_{m})$, $Y_{m}$, $\realN_{+}^{N\times{}K}$]
    \STATE $\muAlgThree = \betaAlgThree / \hat{\nu}_{+} = (\vecOne - \rho)\cdot \betaAlgThree / \hat{\nu}$
    \STATE $\tauAlgThree = Q\betaAlgThree$
    \RETURN $\betaAlgThree, \muAlgThree, \tauAlgThree$
\end{algorithmic}
\end{algorithm}

As with Algorithm~\ref{alg:new_approach_nu_beta_mu_tau}, the extinction cross-section $\beta$ has to be upper bounded so that the loss function (line 4 of Algorithm~\ref{alg:different_approach_beta}) is strictly convex. Let $\mathcal{B}$ be the set over which the loss function is strictly convex. The procedure given in \S\ref{subsubsec:uniqueness_scatter_coeff} to define thet set $\mathcal{M}$, can be used to define the set $\mathcal{B}$ with some minor modifications.

\section{Experimental results}
\label{sec:experimental_results}

\begin{figure}[!ht]
    \centering
    \includegraphics[width=\linewidth]{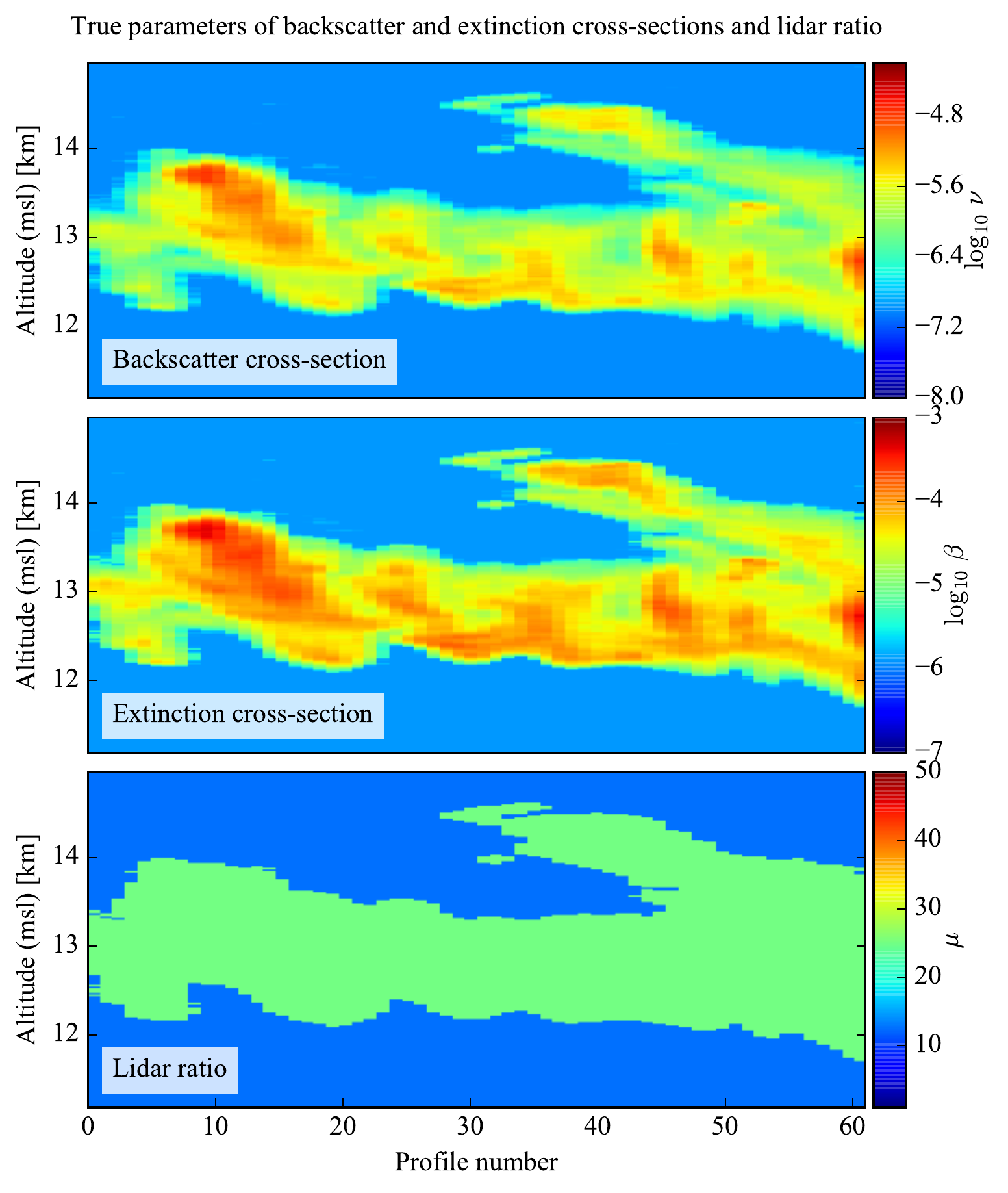}
    \caption{This figure shows the parallel backscatter (top image) and extinction (middle image) cross-sections that were used in the simulation experiments. The bottom image shows the corresponding lidar ratio. The lidar ratio inside the cloud is~25, and for the clear-air the lidar ratio is~40. The linear depolarization coefficients~$\rho$ for the whole scene was set to zero.}
    \label{fig:sim_true_parameters}
\end{figure}

A synthetic dataset was used to juxtapose the performances of Algorithms~\ref{alg:standard_approach_algorithm}, \ref{alg:new_approach_nu_beta_mu_tau} and~\ref{alg:different_approach_beta}. For this section these algorithms will be referred to as the standard (Algorithm~\ref{alg:standard_approach_algorithm}), new (Algorithm~\ref{alg:new_approach_nu_beta_mu_tau}) and alternative (Algorithm~\ref{alg:different_approach_beta}) approaches. The synthetic dataset consists of a cirrus cloud; Figure~\ref{fig:sim_true_parameters} shows the parameters of the synthetic dataset. 
The top and middel images show the parallel backscatter and extinction cross-sections, and the bottom image shows the corresponding lidar ratio. The lidar ratio inside the cloud was set to~25, and outside the cloud it was set to~40. The linear depolarization coefficients $\rho$ for the whole scene was set to zero. Thus, for the simulation experiments the parallel backscatter and backscatter cross-sections are equal to each other.

The synthetic dataset was created using real observations; see~\S\ref{subsec:case_study_two}. The corresponding parallel backscatter cross-section for the simulation scene was generated from an inferred parallel backscatter cross-section using the standard approach, and the residual noise in the inferred parallel backscatter cross-section was suppressed using a Savitzky-Golay (SG) filter that was applied on the range axis.

Two separate experiments were conducted. In the first experiment the parallel backscatter cross-section was inferred, and the performance of the standard and new approaches were compared against each other. And in the second experiment the extinction cross-section was inferred, and the performance of the standard, new and alternative approaches were compared against each other. The Root Mean Square Error (RMSE), which is defined in the following subsection, was used to measure the performance of the different approaches, along with the bias and standard deviation. 

For both experiments the parameters shown in Figure~\ref{fig:sim_true_parameters} were used, although the SNR of the experiments differed. For the first experiment the SNR of a scene represents a real dataset where the spatial resolution is 7.5m by~2.5s. And for the second experiment the SNR was increased by a factor of 48, which represents a scene where the spatial resolution of the original scene is decreased to 7.5m by~120s. This was achieved by oversampling each column of Figure~\ref{fig:sim_true_parameters} by a factor of 48 and the noisy observations of the 48 oversamples columns were accumulated together. The SNR of the second experiment had to be increased so that the extinction cross-section can be inferred, since with the HSRL we work with the molecular-channel observations have much lower SNR values compared to the combined-channel~\cite{eloranta2014high}. The background energies at the native temporal resolution of 2.5s, for both experiments, were set to be 100~times larger compared to a nighttime scene. For the combined- and molecular-channel background energies were 119.29 and 21.46. And since the second experiment represents a scene with a lower spatial resolution, the combined- and molecular-channel background energies were 5725.69 and 1030.18.

The SG filter was used as the lowpass filter for the standard approach, for both the temporal and range axes \cite{eloranta2014high,schafer2010frequency}. The temporal axis SG filter used first order polynomials, and the window size was 9 profiles. The range axis SG filter also used first order polynomials, and the window size was 101 altitude bins. These SG filter parameters were selected based on the minimization of the extinction cross-section RMSE.

\subsection{Error performance measurements}
\label{subsec:error_performance_measurements}

The RMSE is used to compare the performance of the different algorithms. The Mean Square Error (MSE), the square of the RMSE, is defined as 
\begin{align}
    & \text{MSE}\,(\nuAlgOne) = \EX\left[\|\nu - \nuAlgOne\|_{F}^{2}\right] \\
    & = \underbrace{\|\nu - \EX[\nuAlgOne]\|_{F}^{2}}_{\text{bias squared}}
      + \underbrace{\EX\left[\|\nuAlgOne - \EX[\nuAlgOne]\|_{F}^{2}\right]}_{\text{variance}}
\end{align}
We also compare the bias and standard deviation of each algorithm. The norm $\|\cdot\|_{F}^{2}$ is know as the Frobenius norm and it is defined as
\begin{equation}
    \|\nu\|_{F}^{2} = \text{trace}\left(\nu^{\mT}\nu\right).
\end{equation}
The simulation experiments were repeated 25 times in order to obtain values for the RMSE. 

\subsection{Experiment one: Backscatter cross-section estimation}
\label{subsec:experiment_one}

In this experiment the performance of the different approaches were compared against each other when inferring for the  backscatter cross-section at a maximum spatial resolution. The experiment for the standard approach was conducted twice, where for the first time no block averaging was done. And for the second time block averaging was done to reduce the noise variance, where the spatial resolution of the backscattered energy images was reduced to~15m by~5s. 

For this specific experiment the optical depth of the new approach was not calculated using the inferred extinction cross-section, since the extinction cross-section was not estimated in this experiment. The optical depth was calculated using the denoised estimates of the combined- and molecular-channels, and these were plugged into~\eqref{eq:standard_retrieve_beta_p}. This change in how the optical depth was calculated, is specific only to this experiment. 

\begin{figure}[!ht]
    \centering
    \includegraphics[width=\linewidth]{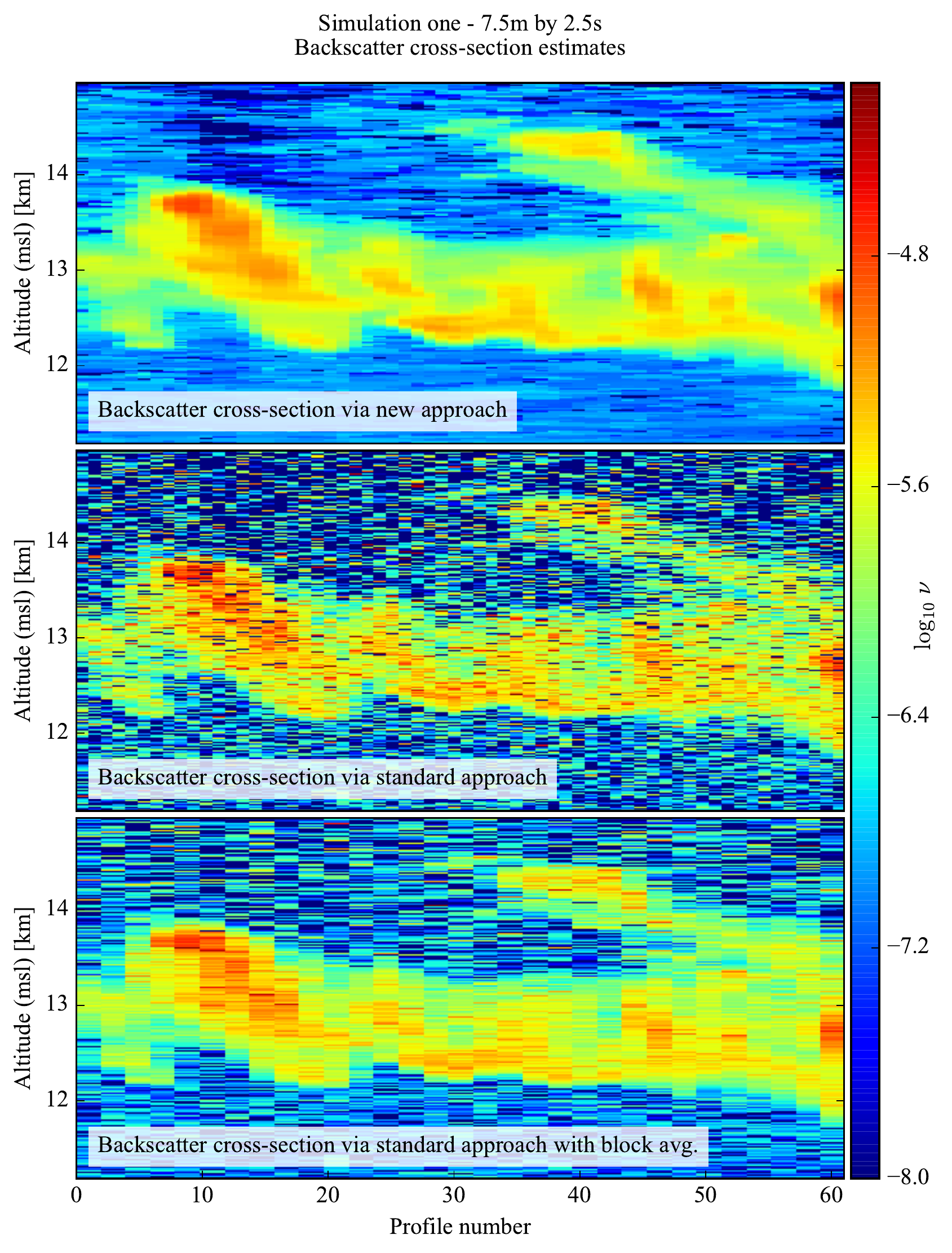}
    \caption{The top image show the inferred backscatter cross-section of the new approach. The middle and bottom images correspond to the standard approach, without and with block averaging. It is clear that the new approach is able to infer the backscatter cross-section with less residual noise (top v.s. middle image) and at a higher spatial resolution (top v.s. bottom image). The new approach is able to preserve the high spatial resolution of the backscatter cross-section due to the Total Variation (TV) smoothness constraint. This smoothness constraint, in conjunction with the noise model, makes the new approach less sensitive to noise.}
    \label{fig:sim1_backscatter_result}
\end{figure}

Table~\ref{tbl:sim_1_results} shows the RMSE values, along with the bias and standard deviation (Std) values, of the different approaches; all the quantities are expressed in the unit of decibel (dB). 
\begin{table}[!ht]
    \caption{This table presents the Root Mean Square Error (RMSE), the bias and standard deviation (Std) results of the first experiment. All the results are expressed in the unit of decibel (dB). From the first four rows we see that the  backscatter cross-section of the new approach attains better performance, compared to the standard approach without and with block averaging. And in regards with the optical depth estimates, the same conclusion follows from the last four rows. This is because new approach use the actual noise model to find an estimates of the parallel backscatter cross-section and optical depth. And it is inferred as an \emph{image}, where the image is constrained to be piecewise constant. In contrast the standard approach infer the parallel backscatter cross-section and optical depth from individual profiles and any spatial and temporal information is not utilized, except for when block averaging is employed which is clearly suboptimal.}
    \label{tbl:sim_1_results}
    \centering
    \begingroup\setlength{\fboxsep}{2pt}
    \colorbox{Gray}{%
    \begin{tabular}{@{}p{3.4cm}lll@{}} \toprule
        \emph{Backscatter cross-section}  &   \emph{RMSE (dB)}  &   \emph{Bias (dB)}  &   \emph{Std (dB)} \\ \midrule
        Standard approach                 &  -8.7740  & -15.7691 & -8.8624 \\
        Standard approach with block avg. &  -25.0532 & -32.0077 & -25.1433 \\      
        New approach                      &  -44.7839 & -45.7323 & -47.0395 \\ 
    \end{tabular}%
    }\endgroup
    
    \begin{tabular}{@{}p{3.4cm}lll@{}} \toprule
        \emph{Optical depth}              &   \emph{RMSE (dB)}  &   \emph{Bias (dB)}  &   \emph{Std (dB)} \\ \midrule
        Standard approach                 &  18.5227 & 11.9207 & 18.4163  \\
        Standard approach with block avg. &  16.0805 & 9.9969  & 15.9445 \\
        New approach                      &  8.5133  & 8.1793  & 4.2833  \\ \arrayrulecolor{color2}\bottomrule
    \end{tabular}
\end{table}
Figure~\ref{fig:sim1_backscatter_result} show results of the inferred backscatter cross-sections of the different approaches. 
The top image show the inferred backscatter cross-section of the new approach. The middle and bottom images correspond to the backscatter cross-section of standard approach without and with block averaging. 

From Table~\ref{tbl:sim_1_results} we see that the new approach is able to infer the backscatter cross-section and optical depth with smaller RMSE values, compared to the standard approach. Furthermore, the corresponding bias and standard deviation of the new approach is less compared to the standard approach. Block averaging did improve the performance of the standard approach, however, the spatial and temporal structure of the cloud is not preserved.

From Figure~\ref{fig:sim1_backscatter_result}, it is clear that the new approach is able to infer the backscatter cross-section with less residual noise (top image v.s. middle image) and at a higher resolution (top image v.s. bottom image). The reason why the new approach is able to achieve these better performance metrics, is because the Total Variation (TV) smoothness constraint is enforced on the recovered backscattered energy and the Poisson noise model is used to measure how good of a ``fit'' the reconstructed backscattered energy is. With the TV smoothness constraint, the underlying backscattered energy is approximated by a spatially piecewise constant function. This allows the estimator of the new approach to preserve the cloud boundaries and the backscatter cross-section variability inside the cloud. 

\subsection{Experiment two: Extinction cross-section estimation}
\label{subsec:experiment_two}

In this experiment the performance of the different algorithms were compared against each other when inferring for the extinction cross-section. In order for standard approach to be able to infer the extinction cross-section accurately for a daytime scene, the image resolution would have been chosen to be much lower than 7.5m by 120s. An image resolution of 7.5m by 120s was chosen, however, to demonstrate what the new approach is able to achieve when compared to the standard approach for a higher resolution image. Also, the simulated scene that was created for this experiment is similar to the real dataset that was used in case study two; see \S\ref{subsec:case_study_two}. We want to use this experiment to show the behaviour of the standard approach when the lidar ratio is inferred. 

Recall that with the new approach the lidar ratio has to be upper bounded so that the loss function is strictly convex; see~\S\ref{subsubsec:uniqueness_scatter_coeff}. If a single upper bound lidar ratio is assumed for the whole scene, the upper bound can be set to be 100. 

\begin{figure}[!ht]
    \centering
    \includegraphics[width=\linewidth]{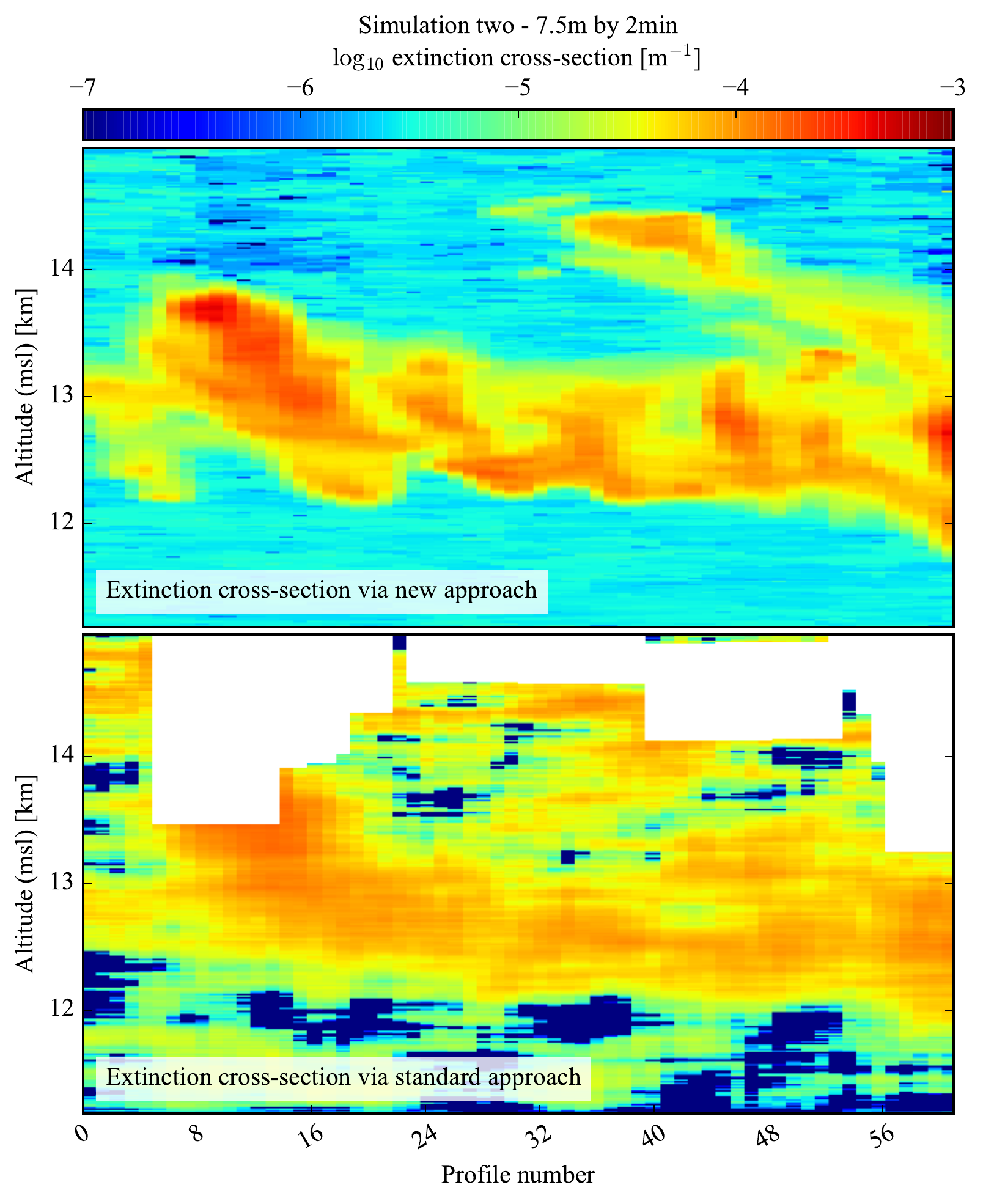}
    \caption{The top image shows the inferred extinction cross-sections of the new approach, and the bottom image shows the extinction cross-section of the standard approach. The white areas in the extinction cross-section image of the standard approach (bottom image), correspond to altitude bins where the inferred transmittance values were negative. And since a lowpass filter is used on the temporal axis, the filter propagated these negative values over multiple profiles. The new approach is able to preserve the high spatial resolution of the extinction cross-section, because 1) it uses the previously estimated backscatter cross-section to infer the lidar ratio and 2) the Total Variation (TV) smoothness constraint is enforced on the lidar ratio. This smoothness constraint, in conjunction with the noise model, makes the new approach resilient against noise.}
    \label{fig:sim2_extinction_result}
\end{figure}

\begin{figure}[!ht]
    \centering
    \includegraphics[width=\linewidth]{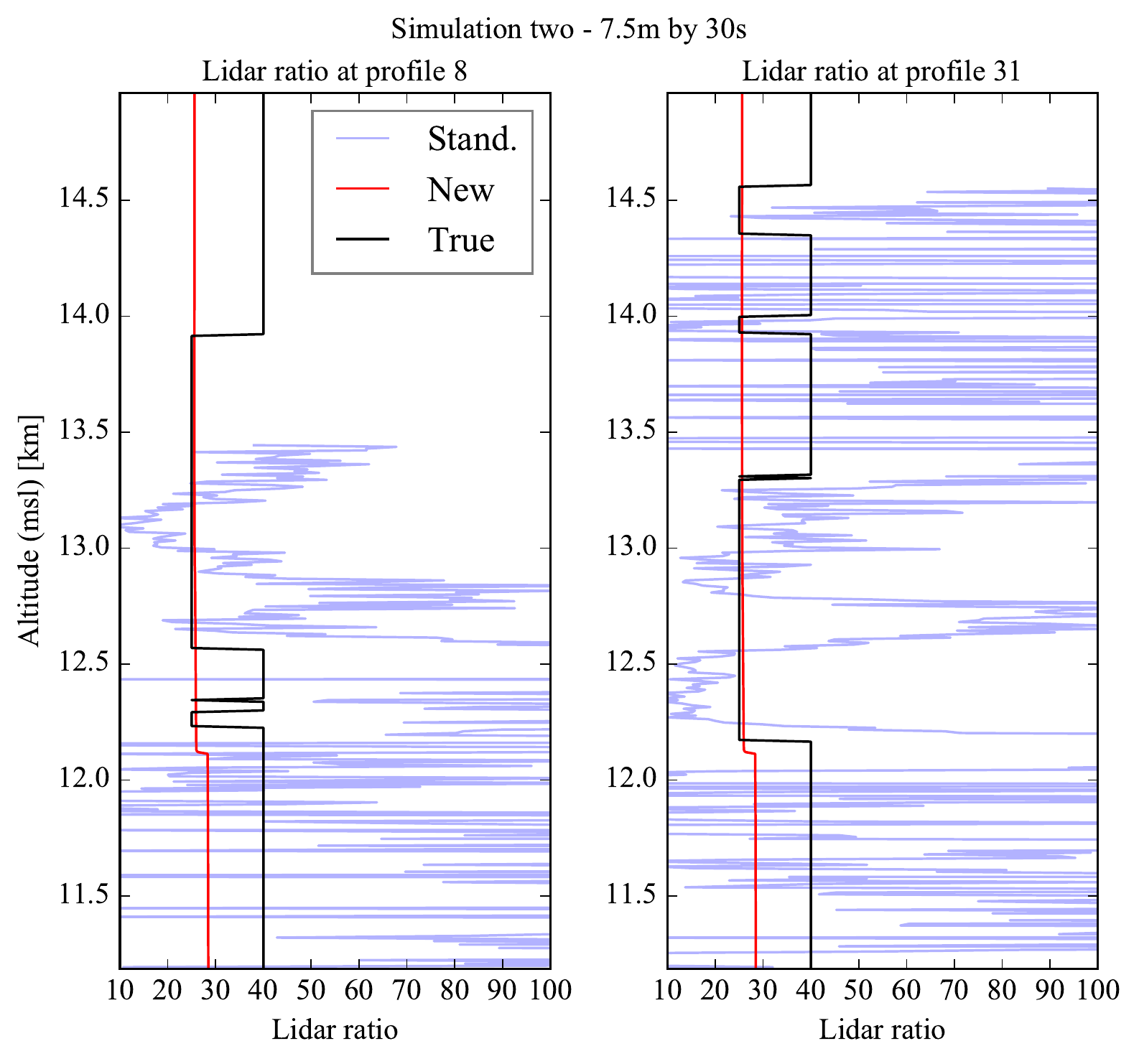}
    \caption{The two images show the inferred lidar ratios of the new (red graphs) and standard (stand., blue graphs) approaches. The black graphs show the true lidar ratios. The lidar ratios beyond~13.5km of profile~8 and beyond~14.5km of profile~31 of the standard approach are missing, because the inferred transmittance values were negative. We see that the new approach is able to infer the lidar ratio to some degree - the accuracy of the clear-air lidar ratio values are off, whereas the standard approach struggles to achieve the same accuracy. Since the Total Variation~(TV) smoothness constraint is enforced on the lidar ratio, the estimates of the lidar ratio are less sensitive to noise. From the rightmost image we see that the standard approach introduces spurious lidar ratios at 12.7km. Such spurious lidar ratios correspond the regions where the extinction cross-section varies; see \S\ref{subsec:spurious_lidar_ratios}.}
    \label{fig:sim2_lidarratio_result}
\end{figure}

Table~\ref{tbl:sim_2_results} shows the RMSE values, along with the bias and standard deviation (Std) values, of the different approaches; all the quantities are expressed in units of dB.
\begin{table}[!ht]
    \caption{This table presents the Root Mean Square Error (RMSE), the bias and standard deviation (Std) results of the second experiment. All the results are expressed in the unit of decibel (dB). From the first three rows, we see that the new approach is able to infer the backscatter cross-section at a higher accuracy compared the standard approach. In all accounts the new approach is able to infer the optical depth, extinction cross-section and lidar ratio with smaller RMSE values compared to the standard and alternative approaches. The new approach is able to achieve these smaller RMSE values for the extinction cross-section, because it uses the previously estimated backscatter cross-section to infer the lidar ratio. And the lidar ratio is constrained to be spatial piecewise constant using the Total Variation (TV) smoothness constraint. It follows then that the inferred optical depth of the new approach also has a higher accuracy, compared to the other approaches.}
    \label{tbl:sim_2_results}
    \centering
    \begingroup\setlength{\fboxsep}{2pt}
    \colorbox{Gray}{%
    \begin{tabular}{@{}p{3.4cm}lll@{}} \toprule
        \emph{Backscatter cross-section} &   \emph{RMSE (dB)}  &   \emph{Bias (dB)}  &   \emph{Std (dB)} \\ \midrule
        Standard approach         &  -33.5874 & -40.3714 & -33.6851 \\
        New approach              &  -44.0111 & -47.5035 & -44.4963 \\
    \end{tabular}%
    }\endgroup
    
    \begin{tabular}{@{}p{3.4cm}lll@{}} \toprule
        \emph{Optical depth}             &   \emph{RMSE (dB)}  &   \emph{Bias (dB)}  &   \emph{Std (dB)} \\ \midrule
        Standard approach         &  14.3185 & 7.9549 & 14.1994 \\
        New approach              &  -3.6500 & -4.8648 & -5.4905 \\     
        Alternative approach      &  2.2543  & 0.5677  & 0.9165 \\
    \end{tabular}
    \begingroup\setlength{\fboxsep}{2pt}
    \colorbox{Gray}{%
    \begin{tabular}{@{}p{3.4cm}lll@{}} \toprule
        \emph{Extinction cross-section} &   \emph{RMSE (dB)}  &   \emph{Bias (dB)}  &   \emph{Std (dB)} \\ \midrule
        Standard approach        & -22.7389 & -23.7649 & -24.8598  \\    
        New approach             & -29.3890 & -33.2679 & -29.7873 \\   
        Alternative approach     & -23.1688 & -24.4253 & -24.9547 \\
    \end{tabular}
    }\endgroup
    
    \begin{tabular}{@{}p{3.4cm}lll@{}} \toprule
        \emph{Lidar ratio}              &   \emph{RMSE (dB)}  &   \emph{Bias (dB)}  &   \emph{Std (dB)} \\ \midrule
        Standard approach        &  74.1082 & 67.1184 & 74.0196 \\
        New approach             &  32.7899 & 32.7673 & 22.8730 \\
        Alternative approach     &  67.6048 & 60.6158 & 67.5161 \\ \arrayrulecolor{color2}\bottomrule
    \end{tabular}
\end{table}
Figure~\ref{fig:sim2_extinction_result} show results of the inferred extinction cross-sections of the new (top image) and~standard (bottom image) approaches.
The white areas in the extinction cross-section image of the standard approach (bottom image), correspond to profiles where the inferred transmittance values were negative. And since a lowpass filter is used on the temporal axis, the filter propagated these negative values over multiple profiles. Figure~\ref{fig:sim2_lidarratio_result} shows the inferred lidar ratios of the new and standard approaches of profile numbers~8 and~32.
The black graphs show the true values of the lidar ratios. The red and blue graphs correspond to the inferred lidar ratios of the new and standard approaches. The lidar ratios beyond~13.5km of profile~8 and beyond~14.5km of profile~31 of the standard approach are missing, because the inferred transmittance values were negative. 

From Table~\ref{tbl:sim_2_results} we see that new approach is able to infer all the unknowns parameters with smaller RMSE values, compared to the standard and alternative approaches. Furthermore, the bias and standard deviation values of the new approach are all smaller compared to those of the other approaches.

From Figure~\ref{fig:sim2_extinction_result} we see that the new approach is able to infer the extinction cross-section in the presence of a large background energy. The standard approach would require more accumulated profiles (i.e. a lower spatial resolution image) to reduce the variance of the noise to a sufficient level until the transmittance values are non-negative. The reason why the new approach does not suffer from the same issue (negative transmittances) as the standard approach, is because it does not try to infer the extinction cross-section with algebraic expressions such as~\eqref{eq:standard_retrieve_beta_p}. The new approach estimates the lidar ratio using the Poisson loss function~\eqref{eq:mu_mle_loss_function} in conjunction with the HSRL model~\eqref{eq:physical_model_molec_channel}, to ``fit'' the estimate on the noisy observations. And by using the HSRL model to do the ``fitting'', the transmittance is always within the bounds of zero and one, as what it should be. 

Figure~\ref{fig:sim2_lidarratio_result} demonstrates that new approach is able to infer the lidar ratios of the cloud to some degree - the accuracy of the clear-air lidar ratio values is off, while the standard approach struggles to achieve the same level of accuracy. The new approach is able to infer these lidar ratios, because it enforces the TV smoothness constraint on the lidar ratio, and it treats the inference problem as a statistical imaging problem; refer to the discussion in the previous section for more information. In constrast, the lowpass filtering of the standard approach is inadequate to preserve the cloud boundaries and the variation of the extinction cross-section inside the cloud. The cloud boundaries and the large variations inside the cloud correspond to high frequency components in the Fourier domain. Hence, the lowpass filters not only suppress the high frequency components, but also contaminate the cloud extinction cross-section values with clear-air extinction cross-section values.

\subsection{How the inferred lidar ratio is influenced by lowpass filters}
\label{subsec:spurious_lidar_ratios}

\begin{figure}[!ht]
    \centering
    \includegraphics[width=\linewidth]{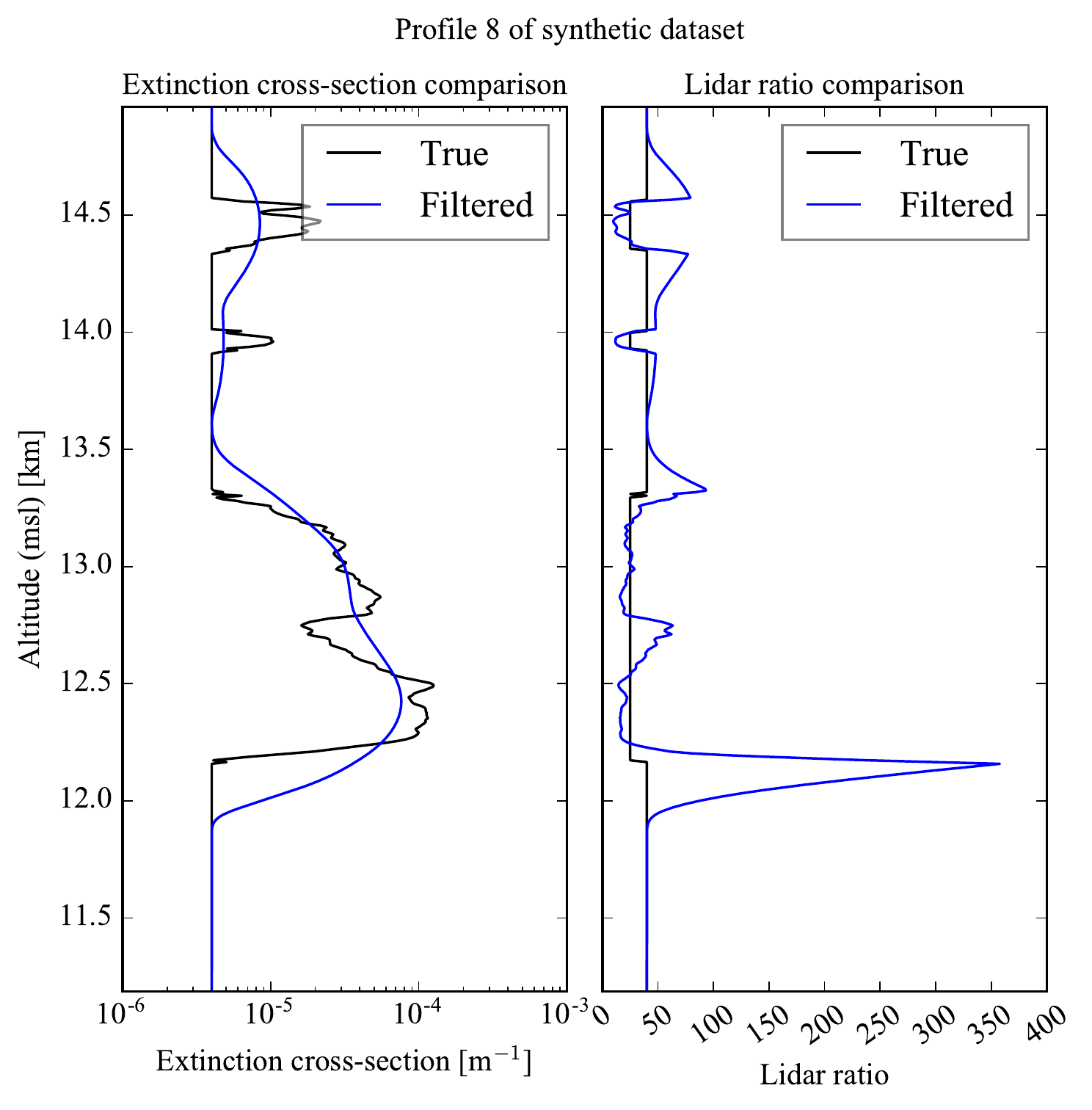}
    \caption{These images explains where the spurious lidar ratios in Figure~\ref{fig:sim2_lidarratio_result} originates from. The leftmost image shows when a lowpass filter is applied on a noiseless optical depth. The spurious lidar ratios on the rightmost image at 12.4km and 12.7km are due to the blurred extinction cross-section values, which are larger than what they should be.}
    \label{fig:lowf_filter_explain}
\end{figure}

From the rightmost image of Figure~\ref{fig:sim2_lidarratio_result} we see that the standard approach introduces spurious lidar ratios at 12.7km. Figure~\ref{fig:lowf_filter_explain} demonstrates that such spurious lidar ratios are most likely induced by the altitude axis lowpass filter, which is used to reduce the residual noise in the inferred optical depth; in this figure noiseless observations were used. 
The altitude axis SG lowpass filter that was used in Figure~\ref{fig:lowf_filter_explain} had a window length of 85 altitude bins. Thus, the bandwidth of the filter is wider than that of Figure~\ref{fig:sim2_extinction_result}~\cite{schafer2010frequency}. The extinction cross-section profile which is used in Figure~\ref{fig:lowf_filter_explain}, corresponds to the lidar ratio profile of the rightmost image of Figure~\ref{fig:sim2_lidarratio_result}. The leftmost image of Figure~\ref{fig:lowf_filter_explain} shows the true extinction cross-section (black graph) and the extinction cross-section that was computed from the lowpass filtered true optical depth (blue graph); see \S\ref{sec:standard_approach}. Notice that at 12.7km the smoothed extinction cross-section is larger than the true extinction cross-section. This implies that the estimated lidar ratio from the smoothed extinction cross-section will be larger than what it is supposed to be, as it is demonstrated in the rightmost image of Figure~\ref{fig:lowf_filter_explain}. The bandwidth of the lowpass filter could be increased to reduce the over-smoothing it induces, but this will be at the expense of not reducing the residual noise to an adequate level.

The lowpass filter smooths out the boundaries of the cirrus clouds at 12.4km of Figure~\ref{fig:lowf_filter_explain}. The implication is that the extinction cross-section of the cloud that increases from $3.5\times{}10^{-6}\,[\meter^{-1}]$ to $5.9\times{}10^{-5}\,[\meter^{-1}]$ over the short altitude interval at 12.4km, is mixed with the clear-air extinction cross-section. Hence, the lidar ratio is unreasonably large at the bottom of the cloud and partially inside the cloud. 

\section{Case studies to evaluate the new approach}
\label{sec:exprmnts_actual_data}

Two case studies of real experimental data are presented to juxtapose the inverted backscatter and extinction cross-sections of Algorithms~\ref{alg:standard_approach_algorithm} and~\ref{alg:new_approach_nu_beta_mu_tau}; these algorithms will be referred to as the standard and new approaches. Both case studies demonstrate that the new approach infer the backscatter and extinction cross-sections at higher spatial and temporal resolutions compared to the standard approach. With both case studies, observations of the Wisconsin-Madison HSRL photon-counting instrument was used, while it was situated at the DOE's Southern Great Plains (SGP) Atmospheric Radiation Measurement (ARM) site during the CHARMS\footnote{Combined HSRL And Raman Measurement Study (CHARMS) was funded by the DOE.} experiment in the summer of~2015. 

We used nighttime scenes to compare the new and standard approaches, since it simplifies the calibration parameters and reduce instrument artifacts. For the second case study, where the extinction cross-section was inferred, the temporal resolution of the molecular-channel noisy image was reduced to increase the SNR of the observations. This is because the HSRL we work with, the molecular-channel observations have much lower SNR values compared to the combined-channel~\cite{eloranta2014high}. The lower SNR is a result of 1) the small molecular backscatter cross-section at the wavelength of the laserpulse which is centered at~532nm, 2) the iodine filter that rejects the particulate backscattered signal (see \S\ref{subsec:hsrl_phy_model}), 3) the telescope diameter of 0.5m, 4) the limit on the laser intensity and 5) the detector quantum inefficiencies. The Wisconsin-Madison HSRL use a Fabry-Perot etalon to reduce the solar background energy, and the telescope diameter is restricted by the diameter of the etalon plates; refer to \cite[\S{}5.3.2]{weitkamp2006lidar} and \cite[Figure~4]{razenkov2010characterization} for more information. The laser intensity is limited due to Federal Aviation Administration (FAA) eye safety standards~\cite[\S{}3.5]{kovalev2004elastic}. The Wisconsin-Madison HSRL use Geiger-Mode Avalanche PhotoDiodes (APD), which have a quantum efficiency of approximately 60\% \cite{razenkov2010characterization}.

\subsection{Case study one - backscatter cross-sections of uniform and non-uniform scenes}
\label{subsec:case_study_one}

\subsubsection{Validation using a uniform HSRL scene}

\begin{figure}[!ht]
\centering
\includegraphics[width=\linewidth]{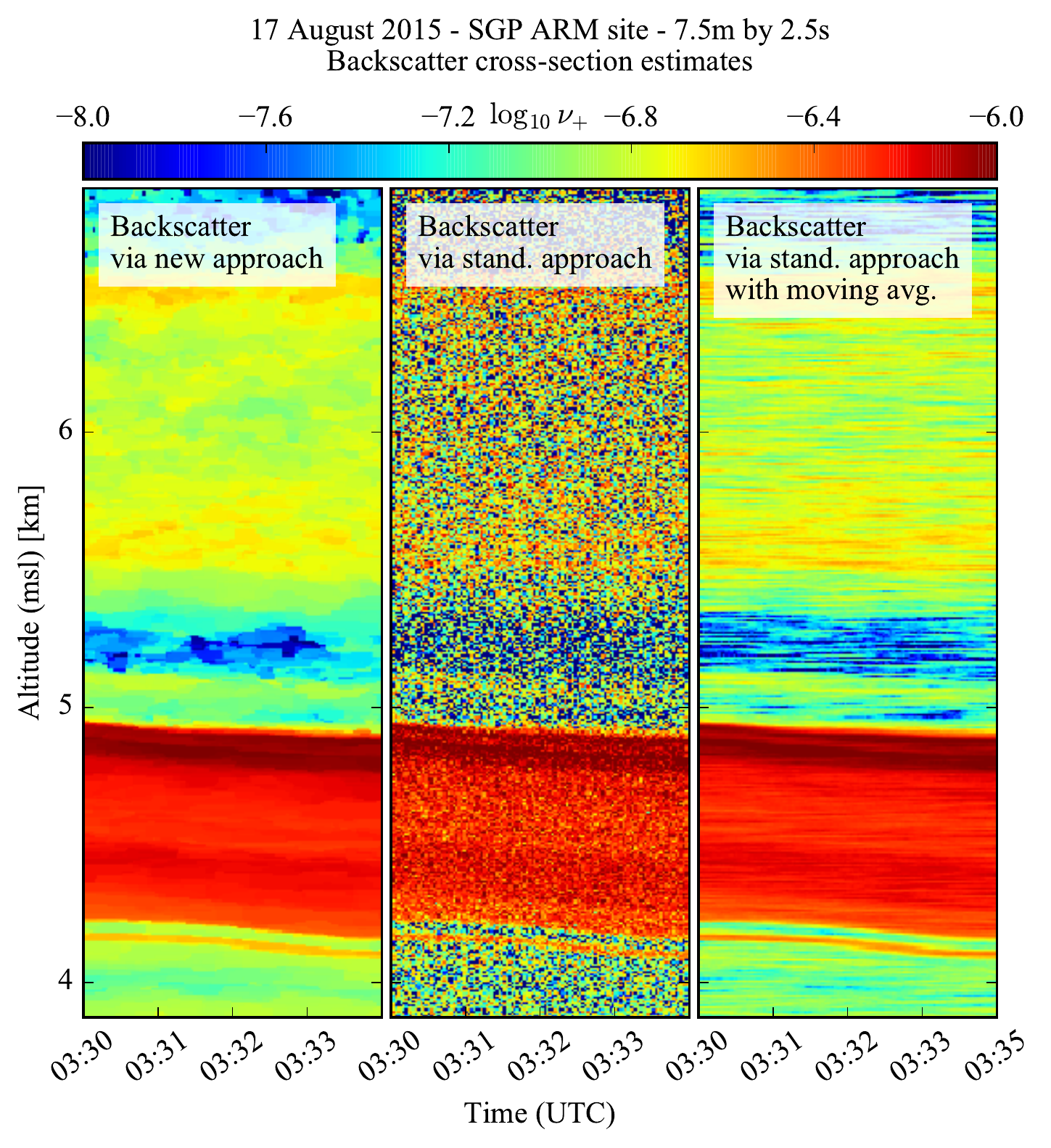}
\caption{In this figure we present a uniform scene, where we expect that the backscatter cross-section of new and standard (stand.) approaches are comparable. It is also used to validate the new approach, in conjunction with Figure~\ref{fig:val1_backscatter_profiles_uniform}. The left image is the backscatter cross-section of the new approach. The middle and the rightmost images show the backscatter cross-section of the standard approach without and with averaging. A moving average filter was used in the rightmost image with a window-width of 1min. This figure demonstrates that the inferred backscatter cross-section of the new approach and standard approach with moving averaging are comparable.}
\label{fig:val1_backscatter_img_uniform}
\end{figure}

\begin{figure}[!ht]
\centering
\includegraphics[width=\linewidth]{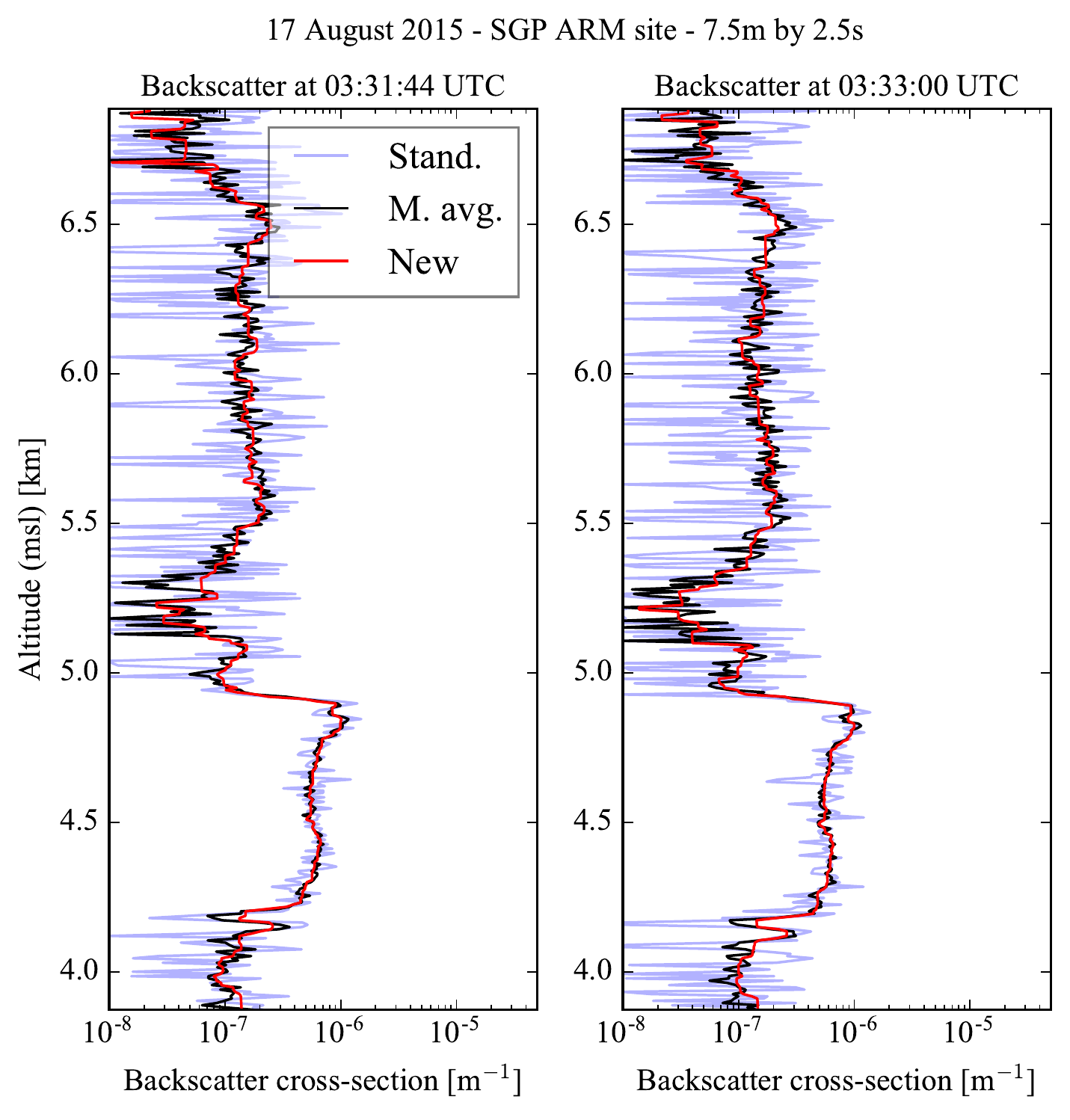}
\caption{The new approach is able to produce an accurate estimate of the backscatter cross-section on real data. The label ``M. avg.'' indicates that a moving average was used to reduce to noise variance, prior to inferring the backscatter cross-section with the standard approach.}
\label{fig:val1_backscatter_profiles_uniform}
\end{figure}

Figures~\ref{fig:val1_backscatter_img_uniform} and~\ref{fig:val1_backscatter_profiles_uniform} validates the  backscatter cross-section of the new approach by using a temporal uniform aerosol scene.
The resolution of the data in both figures is~7.5m by~2.5s. The leftmost image of Figure~\ref{fig:val1_backscatter_img_uniform} shows the inferred backscatter cross-section of the new approach. The middle image is the inverted backscatter of the standard approach with no averaging. In the rightmost image a moving average filter was used with the standard approach, which had a window-width of 1min.

Since the HSRL scene is temporally uniform, it is expected that backscatter cross-section of the standard approach is accurate and the backscatter cross-section of the new approach should be similar to it. To confirm this Figure~\ref{fig:val1_backscatter_profiles_uniform} shows two columns of Figure~\ref{fig:val1_backscatter_img_uniform}. 
For the most part, the backscatter cross-section of the new and standard approaches are comparable. There are specific instances where the new approach potentially over-smooth the backscatter cross-section, for example at~5.5km and~6.5km of the rightmost image.

\subsubsection{Demonstration of preserving the non-uniformity of a cloud}
\label{subsubsec:demonstrate_preserve_boundaries}

\begin{figure}[!ht]
\centering
\includegraphics[width=\linewidth]{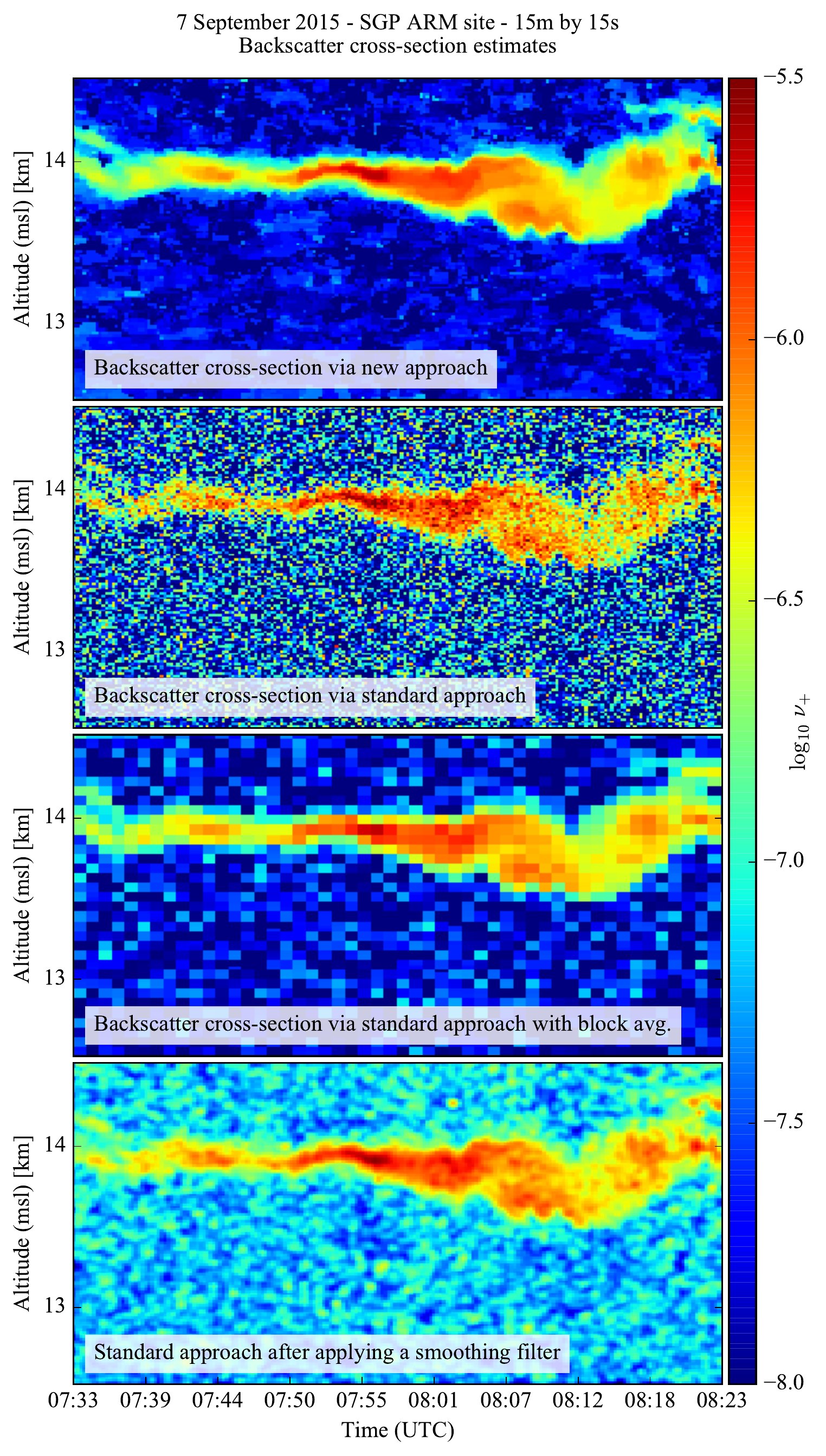}
\caption{The top image shows the backscatter cross-section inverted with the new approach. The new approach is able to preserve the non-uniform spatial structure of a cirrus cloud's backscatter cross-section, while being less sensitive to noise compared to the standard approach. The second image shows the backscatter cross-section inverted with the standard approach at the same resolution as the first image, which is more noisy compared to the first image. The third image shows a block averaged version of the backscatter cross-section of the standard approach, where the resolution is~60m by~60s. Although there is less residual noise in the third image compared to the second image, the spatial structure of the cloud is not well preserved compared to the new approach. The last image shows the inverted backscatter cross-section using the standard approach (second image from the top), after an edge-preserving smoothing filter was applied to suppress the residual noise. The smoothing filter is unable to reduce the noise of the clear-air regions.}
\label{fig:val1_backscatter}
\end{figure}

Figure~\ref{fig:val1_backscatter} demonstrates that the new approach is able to preserve the non-uniform spatial structure of a cirrus cloud's backscatter cross-section at a higher resolution compared to the standard approach, while being less sensitive to noise compared to the standard approach. 
The top image of Figure~\ref{fig:val1_backscatter} shows the backscatter cross-section inverted with the new approach. The second image shows the inverted backscatter cross-section of the standard approach. Both these images have a resolution of~15m by~15s. The third image shows a version of the standard approach backscatter cross-section, where the photon counting observations were block averaged, down to a resolution of 60m by 60s, to reduce the noise variance. The last image shows the inverted backscatter cross-section of the standard approach (second image of Figure~\ref{fig:val1_backscatter}), after an edge-preserving smoothing filter was applied on it to suppress the residual noise; the smoothing filter that was used is a Bilateral filter~\cite{tomasi1998bilateral,scikit-image}. A Bilateral filter modulates the smoothing kernel (i.e. Gaussian function) with a weighting function so that it preserves discontinuities in the region where the filtering is done.

When we compare the top image with the second image, we deduce that the new approach is less sensitive to noise compared to the standard approach. When the top image is juxtaposed against the third image, we deduce that new approach is able to preserve the non-uniform structure at a higher resolution compared to standard approach with block-averaging. The bottom image in comparison with the top image, demonstrates that the smoothing filter is unable to reduce the noise of the clear-air regions.

As it was discussed in \S\ref{subsec:experiment_one}, the reason why the new approach yields better results compared to the standard approach, is because the Poisson noise model is used to ``fit'' the reconstructed backscattered energy on the noisy observations as an image while constraining it be spatially piecewise constant. In other words, each pixel of the backscattered energy of say the combined-channel is not estimated independently, but is estimated in relation to its surrounding backscattered energy pixels and the relation is controlled by the piecewise constant constraint. This allows the new approach to preserve the discontinuities in the backscatter cross-section, while being also less sensitive to the noise compared to the standard approach. Hence, with the new approach the backscatter cross-section can be inferred at a higher resolution than the standard approach. In constrast the standard approach does not take the Poisson noise model in account at all. Furthermore, preprocessing such as block averaging does not accurately maintain the spatial structure of a cirrus cloud. Postprocessing such as lowpass filtering does not necessarily reduce the overall residual noise in a backscatter cross-section image, as it was demonstrated in Figure~\ref{fig:val1_backscatter}.

\subsection{Case study two - extinction cross-section of a non-uniform cirrus cloud}
\label{subsec:case_study_two}

\begin{figure}[!ht]
\centering
\includegraphics[width=\linewidth]{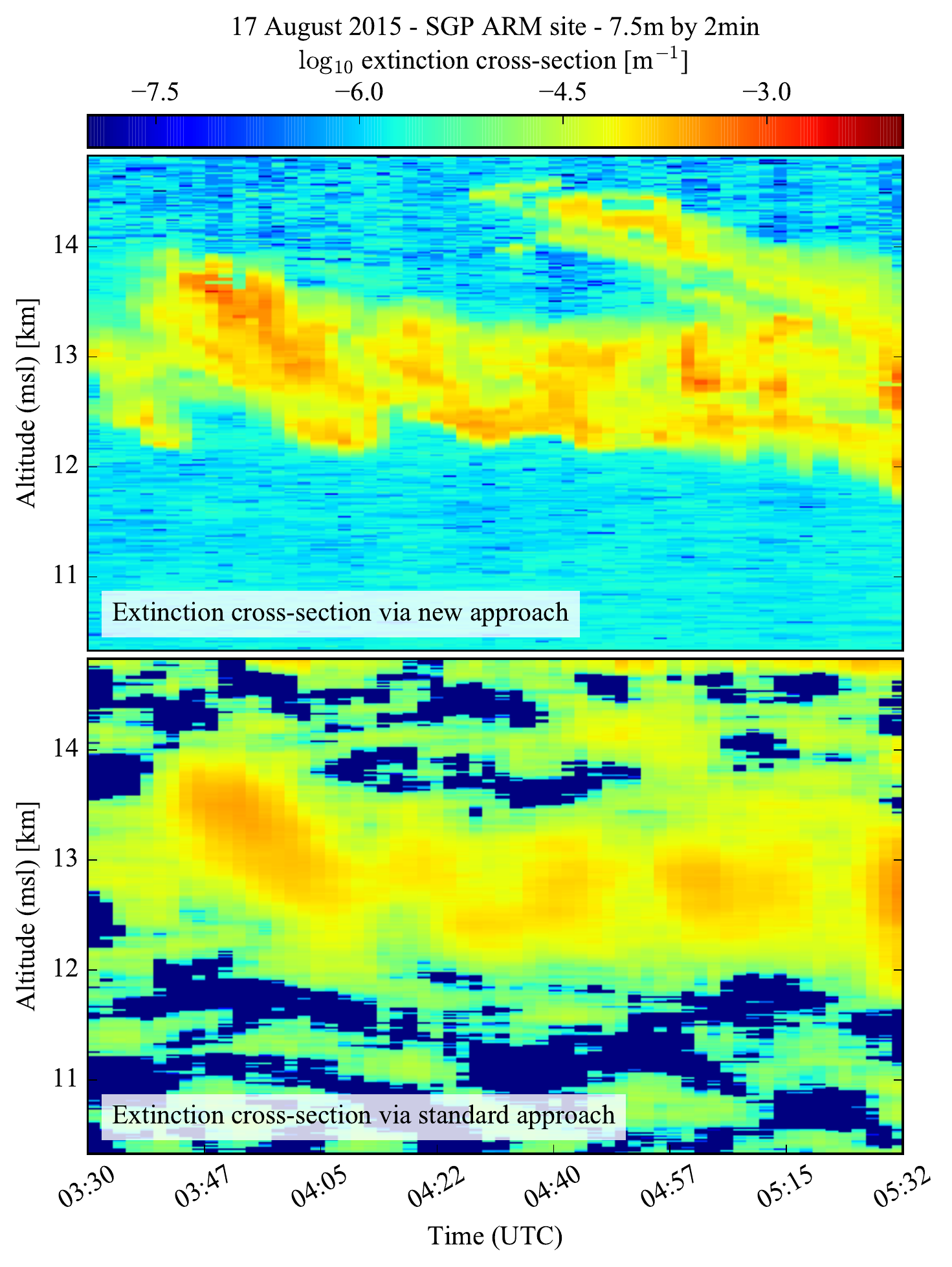}
\caption{The two images show the inverted extinction cross-section of the new and standard approaches, where the topmost image corresponds to the new approach. The new approach is able to preserve the cloud boundaries and the variation of the extinction cross-section inside the cloud, which is due to the Total Variation (TV) smoothness constraint that is enforced on the inferred lidar ratio. Since the standard approach uses lowpass filters to reduce the residual noise of the inferred optical depth on both the temporal and altitudes axes, the cloud boundaries and the extinction cross-section variation inside the cloud are smoothed out.}
\label{fig:val4_scatter_img}
\end{figure}

The inversion of the extinction cross-section and lidar ratio is very sensitive to noise. For the standard approach the extinction is inverted by taking the derivative of the optical depth. The derivative is very sensitive to residual noise in the inferred optical depth. Figure~\ref{fig:val4_scatter_img} presents a case study for a non-uniform cirrus cloud, where we applied the new (top image) and standard (bottom image) approaches over the same photon-counting noisy image.
Notice the new approach is able to resolve the fine scale features of the extinction cross-section, compared to the standard approach. This is because the new approach use the inferred backscatter cross-section to constrain the inversion of the lidar ratio, and the statistical estimator it uses exploit the spatial and temporal correlations in the image. The TV smoothness constraint that the statistical estimator uses, preserves the variation inside the cloud and the cloud boundaries. In contrast the standard approach use lowpass filters to suppress residual noise, resulting in over smoothed cloud boundaries and fine scale features of the extinction cross-section. 

\begin{figure}[!ht]
\centering
\includegraphics[width=\linewidth]{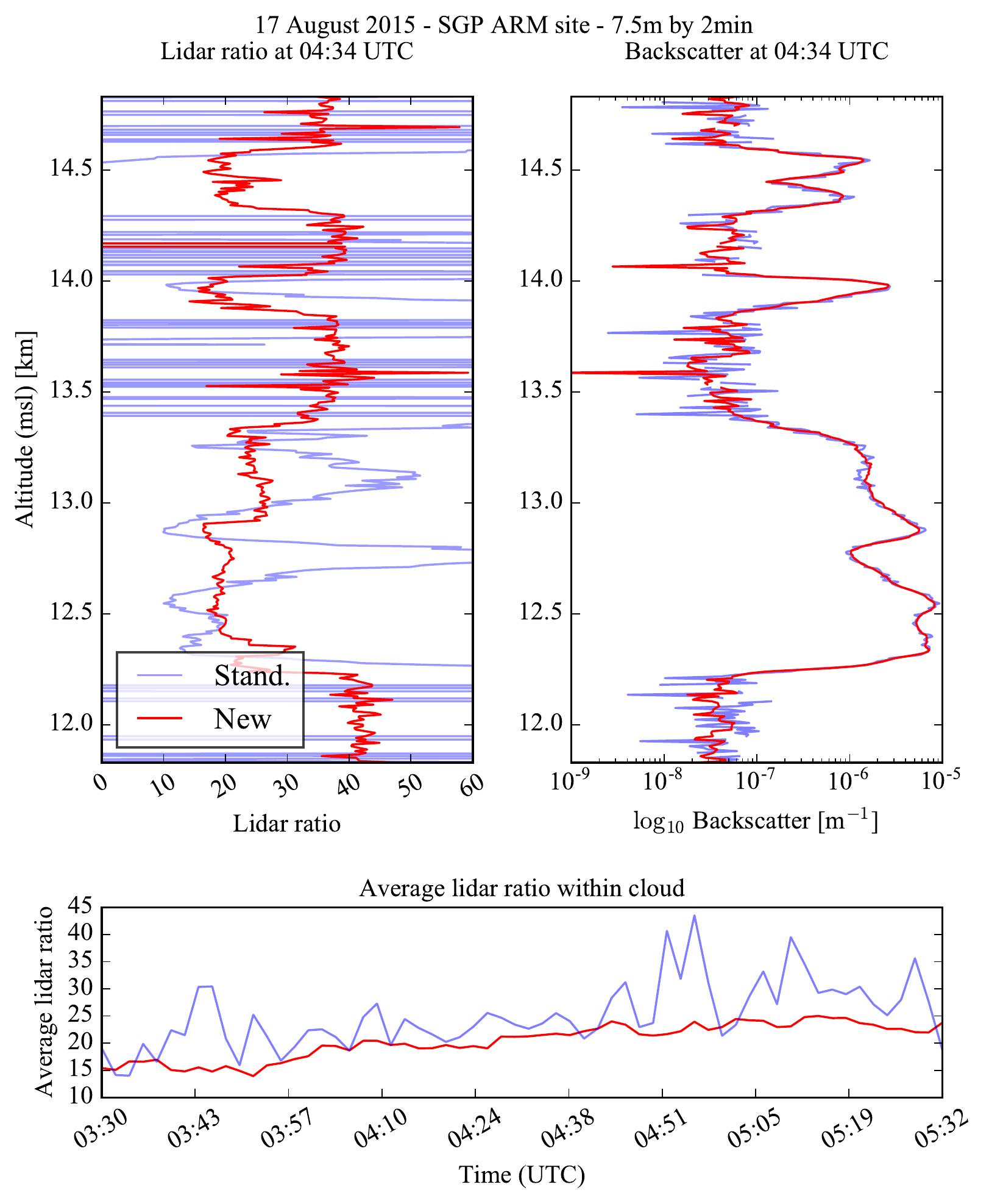} 
\caption{This figure corresponds to Figure~\ref{fig:val4_scatter_img}. The top left image shows that the standard (stand.) approach produces spurious lidar ratio estimates, whereas the lidar ratios of new approach have more realistic variations. As presented in the simulations in \S\ref{subsec:spurious_lidar_ratios}, the spurious lidar ratio estimates are likely a result of the non-uniform structure of the cloud. The change in lidar ratio of the new approach above 13.8km, is likely due to the piecewise constant approximation of the TV-PMLE estimator. The top right image shows the corresponding inverted backscatter cross-sections of the new and standard approaches; these are included to show that the backscatter cross-sections for the two approaches are mostly in agreement. The bottommost image show the average lidar ratios horizontally across the cloud, where the each average was calculated within the cloud. Notice that there is a significant bias between the standard and new approach. As presented in the simulations, in \S\ref{subsec:spurious_lidar_ratios}, this variability and bias are likely a result of susceptibility of the standard approach to the non-uniform structure of the cloud.}
\label{fig:val4_lidarr_woscatter}
\end{figure}

First order polynomials were used by the Savitzky-Golay (SG) filter in Figure~\ref{fig:val4_scatter_img} for the standard approach, where the time and altitude windows were 10min and 637.5m~\cite{eloranta2014high,schafer2010frequency}. Recall that with the new approach the lidar ratio has to be upper bounded so that the loss function is strictly convex; see \S\ref{subsubsec:reparameterization_beta} and \S\ref{subsubsec:uniqueness_scatter_coeff}. If a single upper bound lidar ratio is assumed for the whole scene for this specific case study, the upper bound is approximately 50. For this case study we did set the lidar ratio upper bound to 100. This is reasonable for this case study, since the vertical averaged lidar ratio per profile is less than 100 according to the bottom image of Figure~\ref{fig:val4_lidarr_woscatter}.

The top left image of Figure~\ref{fig:val4_lidarr_woscatter} shows a single column, at 4:34~UTC, of the inferred lidar ratio of the new and standard approaches. 
The top right image is the corresponding backscatter cross-sections of the different approaches; these are included to show that the backscatter cross-sections for the two approaches are mostly in agreement. The bottom image show the vertically averaged lidar ratio of the new and standard approaches, where each average was calculated within the cloud. Notice that there is a significant bias between the standard and new approach, with the standard approach showing more variability and bias relative to the new approach. As presented in the simulations, in \S\ref{subsec:spurious_lidar_ratios}, this variability and bias are likely a result of susceptibility of the standard approach to the non-uniform structure of the cloud. This result demonstrates that not only do we gain the ability the resolve finer features of the extinction cross-section, but the new approach likely produces less biased lidar ratios and extinction cross-sections. 

\section{Conclusion and future work}
\label{sec:conclussion}

We presented a new lidar inversion technique adapted from techniques originally developed for medical imaging. This new approach is able to reduce the random noise while maintaining the spatial and temporal resolutions of the observations, by utilizing the spatial and temporal correlations in images in conjunction with an accurate noise model. By using the data in an entire image in conjunction with the physical noise and lidar models, the new techniques are able to better separate the correlated signal (i.e. clouds and aerosols) from the random noise, resulting from both the detector and solarbackground radiation. The standard approach is to simply average the data in small increments to reduce the noise.

We applied the new and standard approaches on both simulated and real data. Based on these results we found that the new approach is able to accurately maintain the spatial and temporal resolution while significantly reducing the noise, without introducing significant biases due the non-linearity of the lidar equation. This is a significant improvement upon the standard approach, which either has to further average the data to increase the SNR or use lowpass filtering to get comparable noise reduction. The greatest benefit of the new approach relative to the standard approach, is when inverting the extinction cross-section and lidar ratio, which are most susceptible to noise. The standard approach needs additional averaging, which increases the biases due to the non-linearity of the lidar equation.

In the current implementation, the new approach requires Poisson noise statistics (a photon lidar system) and significant computational processing resources. Our future work will focus on improving the computational efficiency of the new technique and apply it to a large ensemble of observations. The next step is to adapt the new techniques to standard backscatter photon lidar instruments, such as the NASA Cloud-Aerosol Transport System (CATS), Cloud Physics Lidar (CPL) and Micro-Pulsed Lidar (MPL)~\cite{mcgill2012cats,mcgill2002cloud,welton2001global}. A much more significant challenge is to adapt the new technique for analog-mode lidar instruments such as CALIPSO, where the noise is non-Poisson~\cite{liu2002simulation}. Once the techniques have been adapted for standard backscatter and analog-mode lidar instruments, the techniques have the potential to greatly improve the scientific utility of systems such as the NASA CATS and NASA-CNES Cloud-Aerosol Lidar and Infrared Pathfinder Satellite Observations~\cite{mcgill2012cats,winker2009overview}, which have lower SNR observations compared to ground based lidar systems.

\appendix 

\begin{center}
  {\bf APPENDIX}
\end{center}

\section{Thinning a Poisson distributed matrix in creating the training and testing matrices}
\label{app:poisson_thinning}

Without loss of generality, suppose $Y$ be a Poisson random variable where each entry is independent from each other
and $\EX_{Y}[Y] = \lambda$. And for a given $Y$, let $Z$ be a Binomial random variable, where $\EX_{Z|Y}[Z | Y = y] =
yp$ and $p\in(0,1)$. The random variables $Y$ and $Z$ are independent from each other. It can be proved that $Z$ is
Poisson distributed with $\EX_{Z,Y}[Z] = \EX_{Y}[\EX_{Z|Y}[Z | Y]] = p\lambda$, where the expectation is taken over the
PMFs of both $Y$ and $Z$.

Now suppose $Y_{\iota}$ is Poisson distributed random matrix. Using the procedure described in the previous paragraph,
we can create a ``training'' matrix $\trnYi$ where $\EX[\trnYi] = p\EX[Y_{\iota}]$. The ``testing'' matrix $\tstYi$ is
created by the removing $\trnYi$ from $Y_{\iota}$. This is achieved by subtracting the two matrices.

\section*{Acknowledgments}

The authors would like to thank Zachary Harmany and Albert Oh for sharing their most recent source code of \cite{oh2013logTV}. We thank Rick Wagener for maintaining the AERONET instruments and making the data available at the SGP ARM site. We also thank Ilya Razenkov, Joe Garcia and Martin Lawson for their dedicated hard work in maintaining and developing the Wisconsin-Madison HSRL system. We also thank the reviewers for their helpful feedback on the paper. 

\bibliographystyle{IEEEtran}
\bibliography{references.bib}

\begin{thebibliography}{10}
\providecommand{\url}[1]{#1}
\csname url@samestyle\endcsname
\providecommand{\newblock}{\relax}
\providecommand{\bibinfo}[2]{#2}
\providecommand{\BIBentrySTDinterwordspacing}{\spaceskip=0pt\relax}
\providecommand{\BIBentryALTinterwordstretchfactor}{4}
\providecommand{\BIBentryALTinterwordspacing}{\spaceskip=\fontdimen2\font plus
\BIBentryALTinterwordstretchfactor\fontdimen3\font minus
  \fontdimen4\font\relax}
\providecommand{\BIBforeignlanguage}[2]{{%
\expandafter\ifx\csname l@#1\endcsname\relax
\typeout{** WARNING: IEEEtran.bst: No hyphenation pattern has been}%
\typeout{** loaded for the language `#1'. Using the pattern for}%
\typeout{** the default language instead.}%
\else
\language=\csname l@#1\endcsname
\fi
#2}}
\providecommand{\BIBdecl}{\relax}
\BIBdecl

\bibitem{winker2009overview}
D.~M. Winker, M.~A. Vaughan, A.~Omar, Y.~Hu, K.~A. Powell, Z.~Liu, W.~H. Hunt,
  and S.~A. Young, ``Overview of the {CALIPSO} mission and {CALIOP} data
  processing algorithms,'' \emph{Journal of Atmospheric and Oceanic
  Technology}, vol.~26, no.~11, pp. 2310--2323, 2009.

\bibitem{hair2008airborne}
J.~W. Hair, C.~A. Hostetler, A.~L. Cook, D.~B. Harper, R.~A. Ferrare, T.~L.
  Mack, W.~Welch, L.~R. Izquierdo, and F.~E. Hovis, ``Airborne high spectral
  resolution lidar for profiling aerosol optical properties,'' \emph{Applied
  Optics}, vol.~47, no.~36, pp. 6734--6752, 2008.

\bibitem{durand2007esa}
Y.~Durand, A.~H{\'e}li{\`e}re, J.~B{\'e}zy, and R.~Meynart, ``The {ESA}
  {E}arth{CARE} mission: results of the {ATLID} instrument pre-developments,''
  in \emph{Proc. SPIE}, vol. 6750, 2007, p. 675015.

\bibitem{mcgill2002cloud}
M.~McGill, D.~Hlavka, W.~Hart, V.~S. Scott, J.~Spinhirne, and B.~Schmid,
  ``Cloud physics lidar: Instrument description and initial measurement
  results,'' \emph{Applied Optics}, vol.~41, no.~18, pp. 3725--3734, 2002.

\bibitem{welton2001global}
E.~J. Welton, J.~R. Campbell, J.~D. Spinhirne, and V.~S. Scott~III, ``Global
  monitoring of clouds and aerosols using a network of micropulse lidar
  systems,'' in \emph{Second International Asia-Pacific Symposium on Remote
  Sensing of the Atmosphere, Environment, and Space}.\hskip 1em plus 0.5em
  minus 0.4em\relax International Society for Optics and Photonics, 2001, pp.
  151--158.

\bibitem{liu2006estimating}
Z.~Liu, W.~Hunt, M.~Vaughan, C.~Hostetler, M.~McGill, K.~Powell, D.~Winker, and
  Y.~Hu, ``Estimating random errors due to shot noise in backscatter lidar
  observations,'' \emph{Applied optics}, vol.~45, no.~18, pp. 4437--4447, 2006.

\bibitem{liu2002simulation}
Z.~Liu and N.~Sugimoto, ``Simulation study for cloud detection with space
  lidars by use of analog detection photomultiplier tubes,'' \emph{Applied
  optics}, vol.~41, no.~9, pp. 1750--1759, 2002.

\bibitem{hastie2009elements}
T.~Hastie, R.~Tibshirani, and J.~Friedman, \emph{The Elements of Statistical
  Learning: Data Mining, Inference, and Prediction}, 2nd~ed., ser. Springer
  Series in Statistics.\hskip 1em plus 0.5em minus 0.4em\relax Springer New
  York, 2009.

\bibitem{rodgers2000inverse}
C.~D. Rodgers, \emph{Inverse methods for atmospheric sounding: theory and
  practice}.\hskip 1em plus 0.5em minus 0.4em\relax World scientific Singapore,
  2000, vol.~2.

\bibitem{eloranta2014high}
E.~Eloranta, ``{H}igh {S}pectral {R}esolution {L}idar {M}easurements of
  {A}tmospheric {E}xtinction: {P}rogress and {C}hallenges,'' in \emph{Aerospace
  Conference, 2014 IEEE}.\hskip 1em plus 0.5em minus 0.4em\relax IEEE, 2014,
  pp. 1--6.

\bibitem{nishizawa2008algorithm}
T.~Nishizawa, N.~Sugimoto, I.~Matsui, A.~Shimizu, B.~Tatarov, and H.~Okamoto,
  ``Algorithm to retrieve aerosol optical properties from
  high-spectral-resolution lidar and polarization mie-scattering lidar
  measurements,'' \emph{Geoscience and Remote Sensing, IEEE Transactions on},
  vol.~46, no.~12, pp. 4094--4103, 2008.

\bibitem{young2008caliopATBDp4}
S.~A. Young, D.~Winker, M.~Vaughan, Y.~Hu, and R.~Kuehn, ``{CALIOP} {A}lgorithm
  {T}heoretical {B}asis {D}ocument, {P}art 4: {E}xtinction retrieval
  algorithms,'' NASA Langley Research Center, Tech. Rep., 2008.

\bibitem{ansmann2007particle}
A.~Ansmann, U.~Wandinger, O.~Le~Rille, D.~Lajas, and A.~G. Straume, ``Particle
  backscatter and extinction profiling with the spaceborne
  high-spectral-resolution doppler lidar aladin: methodology and simulations,''
  \emph{Applied optics}, vol.~46, no.~26, pp. 6606--6622, 2007.

\bibitem{oh2013logTV}
A.~Oh, Z.~Harmany, and R.~Willett, ``Logarithmic total variation regularization
  for cross-validation in photon-limited imaging,'' in \emph{Image Processing
  (ICIP), 2013 20th IEEE International Conference on}, Sept 2013, pp. 484--488.

\bibitem{harmany2012spiral}
Z.~T. Harmany, R.~F. Marcia, and R.~M. Willett, ``This is {SPIRAL}-{TAP}:
  {S}parse {P}oisson {I}ntensity {R}econstruction {AL}gorithms---{T}heory and
  {P}ractice,'' \emph{Image Processing, IEEE Transactions on}, vol.~21, no.~3,
  pp. 1084--1096, 2012.

\bibitem{willett2003platelets}
R.~M. Willett and R.~D. Nowak, ``Platelets: a multiscale approach for
  recovering edges and surfaces in photon-limited medical imaging,''
  \emph{Medical Imaging, IEEE Transactions on}, vol.~22, no.~3, pp. 332--350,
  2003.

\bibitem{weitkamp2006lidar}
C.~Weitkamp, \emph{Lidar: range-resolved optical remote sensing of the
  atmosphere}.\hskip 1em plus 0.5em minus 0.4em\relax Springer Science \&
  Business, 2006, vol. 102.

\bibitem{petty2006first}
G.~W. Petty, \emph{A first course in atmospheric radiation}.\hskip 1em plus
  0.5em minus 0.4em\relax Sundog Pub, 2006.

\bibitem{eloranta2005high}
E.~E. Eloranta, \emph{High spectral resolution lidar}.\hskip 1em plus 0.5em
  minus 0.4em\relax Springer, 2005.

\bibitem{hsrldev-www}
\BIBentryALTinterwordspacing
E.~Eloranta. (2014, 10) http://hsrl.ssec.wisc.edu. [Online]. Available:
  \url{http://hsrl.ssec.wisc.edu}
\BIBentrySTDinterwordspacing

\bibitem{razenkov2010characterization}
I.~Razenkov, ``Characterization of a geiger-mode avalanche photodiode detector
  for high spectral resolution lidar,'' Ph.D. dissertation, UNIVERSITY OF
  WISCONSIN-MADISON, 2010.

\bibitem{willett2007multiscale}
R.~M. Willett and R.~D. Nowak, ``Multiscale poisson intensity and density
  estimation,'' \emph{Information Theory, IEEE Transactions on}, vol.~53,
  no.~9, pp. 3171--3187, 2007.

\bibitem{eggermont2001maximum}
P.~P.~B. Eggermont, V.~N. LaRiccia, and V.~LaRiccia, \emph{Maximum penalized
  likelihood estimation}.\hskip 1em plus 0.5em minus 0.4em\relax Springer,
  2001, vol.~1.

\bibitem{schafer2010frequency}
R.~W. Schafer, ``{O}n the {F}requency-{D}omain {P}roperties of
  {S}avitzky-{G}olay {F}ilters,'' in \emph{Proc. 2011 DSP/SPE Workshop}, 2010,
  pp. 54--59.

\bibitem{pornsawad2012retrieval}
P.~Pornsawad, G.~D'Amico, C.~B{\"o}ckmann, A.~Amodeo, and G.~Pappalardo,
  ``Retrieval of aerosol extinction coefficient profiles from raman lidar data
  by inversion method,'' \emph{Applied optics}, vol.~51, no.~12, pp.
  2035--2044, 2012.

\bibitem{hansen1999curve}
P.~C. Hansen, \emph{The L-curve and its use in the numerical treatment of
  inverse problems}.\hskip 1em plus 0.5em minus 0.4em\relax IMM, Department of
  Mathematical Modelling, Technical Universityof Denmark, 1999.

\bibitem{tian2014improved}
P.~Tian, X.~Cao, J.~Liang, L.~Zhang, N.~Yi, L.~Wang, and X.~Cheng, ``Improved
  empirical mode decomposition based denoising method for lidar signals,''
  \emph{Optics Communications}, vol. 325, pp. 54--59, 2014.

\bibitem{cimel2015ce370-2}
\BIBentryALTinterwordspacing
Cimel. (2015, 6) Caml : Cloud and aerosol micro lidar. [Online]. Available:
  \url{http://support.cimel.fr/photo/pdf/ce370_us.pdf}
\BIBentrySTDinterwordspacing

\bibitem{o1986statistical}
F.~O'Sullivan, ``A statistical perspective on ill-posed inverse problems,''
  \emph{Statistical science}, pp. 502--518, 1986.

\bibitem{needell2013stable}
D.~Needell and R.~Ward, ``Stable image reconstruction using total variation
  minimization,'' \emph{SIAM Journal on Imaging Sciences}, vol.~6, no.~2, pp.
  1035--1058, 2013.

\bibitem{chan2005image}
T.~F. Chan and J.~J. Shen, \emph{Image processing and analysis: variational,
  PDE, wavelet, and stochastic methods}.\hskip 1em plus 0.5em minus 0.4em\relax
  Siam, 2005.

\bibitem{kovalev2004elastic}
\BIBentryALTinterwordspacing
V.~Kovalev and W.~Eichinger, \emph{Elastic Lidar: Theory, Practice, and
  Analysis Methods}.\hskip 1em plus 0.5em minus 0.4em\relax Wiley, 2004.
  [Online]. Available: \url{https://books.google.co.kr/books?id=C17mLdkhXD8C}
\BIBentrySTDinterwordspacing

\bibitem{tomasi1998bilateral}
C.~Tomasi and R.~Manduchi, ``Bilateral filtering for gray and color images,''
  in \emph{Computer Vision, 1998. Sixth International Conference on}.\hskip 1em
  plus 0.5em minus 0.4em\relax IEEE, 1998, pp. 839--846.

\bibitem{scikit-image}
\BIBentryALTinterwordspacing
S.~van~der Walt, J.~L. {S}ch\"onberger, J.~{Nunez-Iglesias}, F.~{B}oulogne,
  J.~D. {W}arner, N.~{Y}ager, E.~{G}ouillart, T.~{Y}u, and the scikit-image
  contributors, ``scikit-image: image processing in {P}ython,'' \emph{PeerJ},
  vol.~2, p. e453, 6 2014. [Online]. Available:
  \url{http://dx.doi.org/10.7717/peerj.453}
\BIBentrySTDinterwordspacing

\bibitem{mcgill2012cats}
M.~McGill, E.~Welton, J.~Yorks, and V.~S. Scott, ``{CATS}: A new earth science
  capability,'' Newsletter, May June 2012.

\end{thebibliography}
\end{document}